\documentclass{article}

\RequirePackage[OT1]{fontenc}
\usepackage{graphicx,amsmath,amsfonts,amssymb,amsthm,psfrag, color}
\usepackage[affil-it]{authblk}
\usepackage{bbm,wrapfig,float}
\usepackage{fullpage}

\newcommand{\ds}{\displaystyle}
\newcommand{\ind}[1]{\mathbbm{1}_{#1}}\newcommand{\dint}{\mathrm{d}}

\newtheorem{thm}{Theorem}[section]

\newtheorem{proposition}[thm]{Proposition}

\newtheorem{rem}[thm]{Remark}
\numberwithin{equation}{section}

\begin{document}
\title{ Simulation of hitting times for Bessel processes with non integer dimension}
\author{Madalina Deaconu 
\thanks{Electronic address: \texttt{Madalina.Deaconu@inria.fr}; Corresponding author}}
\affil{Inria, Villers-l{\`e}s-Nancy, F-54600, France;\\
Universit{\'e} de Lorraine, CNRS, Institut Elie Cartan de Lorraine - UMR 7502,\\
Vandoeuvre-l{\`e}s-Nancy, F-54506, France; }
\author{Samuel Herrmann \thanks{Electronic address: \texttt{Samuel.Herrmann@u-bourgogne.fr}; supported by Conseil Regional de Bourgogne (contract no. 2012-9201AAO047S01283)}}
\affil{Institut de Math{\'e}matiques de Bourgogne (IMB) - UMR 5584,\\
   Universit{\'e} de Bourgogne, B.P. 47 870,\\
   21078 Dijon Cedex, France }

\maketitle

\begin{abstract}
In this paper we pursue and complete the study of the simulation of the hitting time of some given boundaries for Bessel processes. These problems are of great interest in many application fields as finance and neurosciences. In a previous work \cite{Deaconu-Herrmann}, the authors introduced a new method for the simulation of hitting times for Bessel processes with integer dimension. The method was based mainly on the explicit formula for the distribution of the hitting time and on the connexion between the Bessel process and the Euclidean norm of the Brownian motion. This method does not apply anymore for a non integer dimension. 
In this paper we consider the simulation of the hitting time of Bessel processes with non integer dimension and provide a new algorithm by using the additivity property of the laws of squared Bessel processes. We split each simulation step in two parts: one is using the integer dimension case and the other one considers hitting time of a Bessel process starting from zero.
\end{abstract}


\section{Introduction}
\par\noindent\null $\,$ \\ 

The aim of this paper is to construct new methods for approaching the hitting time of a given level for Bessel processes of one non-integer dimension. Diffusion hitting times are important quantities in many fields of sciences and applications, such as mathematical science, finance, geophysics or neurosciences. A typical example is the study of path dependent exotic options as barrier options in finance.

On one hand, analytic expressions for hitting time densities are well known and studied only in some very particular situations.
On the other hand, the study of the approximation of the hitting times for stochastic differential equations is an active area of research since very few results exist already. For the Brownian motion, we can approach this quantity simply by using Gaussian random variables. Four alternatives for dealing with the characterization of hitting times in the Brownian case exist : Monte Carlo and Euler based methods, Volterra integral equation for Gaussian Markov processes, series expansion as performed by Durbin \cite{durbin85} or partial differential equation approaches which are based on the explicit form of the probability distribution function of the Brownian motion. These methods do not apply in the general diffusion case as they rely on this explicit form.
For the general diffusion case, very few studies in this direction exist. The only methods that can be used are the Monte Carlo method and time splitting method like the Euler scheme. Some works have already been done in the context of smooth drift and diffusion coefficient by Gobet and Menozzi \cite{ gobet2001,
gobet-menozzi2010}. A recent paper by Ichiba and Kardaras \cite{ichiba-kardaras2011} also uses the representation of the passage density as the mean of a three dimensional Brownian motion.

Our study focuses on the numerical approach for the hitting time of a Bessel process. The main point of this work is the introduction of an iterative procedure for approaching the hitting time by using the structure and the particular properties of the Bessel process. In particular, our approach avoids splitting time methods.
More precisely, we consider the simulation for the hitting time of Bessel processes with non integer dimension and construct a new algorithm by using the additivity property of the laws of squared Bessel processes. We split each simulation step in two parts~: one uses the integer dimension case and the other considers the hitting time of a Bessel process starting from zero. %

The paper is organized as follows. Sections 2, 3 and 4 introduce some generalities about Bessel processes and the properties needed in the paper. Section 5
gives a description of the hitting time of a particular boundary for a Bessel process starting from 0. In section 6, the algorithm is provided and the convergence discussed. Finally, section 7 gives some numerical results. \\ 
\section{Research topic - Bessel processes}
\par\noindent\null $\,$ \\ 

We consider the $\delta$-dimensional Bessel process ($\delta>1$), starting from $x_0$, which is the solution of the following stochastic differential equation
\begin{equation}
\label{bessel-delta}
\left\{
\begin{array}{ll}
X^{\delta,x_0}_t& = X^{\delta,x_0}_0+\ds\frac{\delta-1}{2}\int_0^t (X^{\delta,x_0}_s)^{-1} \dint s + B_t,\\
X^{\delta,x_0}_0&=x_0, \qquad x_0\geq0.\\
\end{array}
\right.
\end{equation}
Let us recall that the Bessel process is characterized either by its dimension $\delta$ or, alternatively, by its index $\nu$ given by $\nu=\delta/2-1$.\\
For a fixed $L > 0$ let us denote by 
\begin{equation}
\label{hitting_time_L}
\tau_L=\inf\{ t\geq 0; X^{\delta,x_0}_t=L \},
\end{equation}
the first time that the process hits the level $L$. In the Bessel process case an explicit form of the Laplace transform of $\tau_L$ exists :
\begin{equation*}
 \mathbb{E}_{x_0}\Big[ e^{-\lambda \tau_L} \Big]=\frac{(x_0)^{-\nu}}{L^{-\nu}}\frac{I_\nu(x_0\sqrt{2\lambda})}{I_\nu(L\sqrt{2\lambda})},\ y>0,
\end{equation*}
here $I_\nu(x)$ denotes the modified Bessel function.
Ciesielsky and Taylor \cite{Ciesielski-Taylor} also proved that for $\delta \in \mathbb{N}$ the tail distribution is given by, when starting from 0  
\begin{equation*}
 \mathbb{P}(\tau_L>t)=\frac{1}{2^{\nu-1}\Gamma(\nu+1)}\,\sum_{k=1}^{\infty}\frac{j_{\nu,k}^{\nu-1}}{\mathcal{J}_{\nu+1}(j_{\nu,k})}\,e^{-\frac{j_{\nu,k}^2}{2L^2}t}.
\end{equation*}
where $\mathcal{J}_\cdot$  is the Bessel function of the first kind, and $j_{\cdot,k}$ is the associated sequence of its positive zeros.
These formulas are restricted to the integer dimension case (see \cite{Hamana} for non integer dimensions) and are obviously miss-adapted and not suited for numerical approaches.

Let us also recall some properties of Bessel processes with respect to their dimension as follows from Revuz and Yor \cite{revuz-yor} or Jeanblanc, Yor and Chesney {\cite{jeanblanc-yor-chesney}.
\begin{enumerate}
\item For $\delta > 2 $ the process  $BESQ^\delta$ is transient and, for $\delta\leq 2$, it is recurrent.
\item For $\delta \geq 2$ the point 0 is polar for $BESQ^\delta$, and for $\delta \leq 1$ it is reached almost surely. 
\item For $\delta =0 $ the point 0 is absorbing.
\item For $0\ < \delta < 2$ the point 0 is instantaneously reflecting.
\end{enumerate}

First, we introduce some relations connecting Bessel processes of different dimensions. The first relation is based on Girsanov's transformation and the second one is a decomposition of the squared Bessel process as a sum of two independent squared Bessel processes.

\section{Absolute continuity}
\par\noindent\null $\,$ \\ 

On the canonical space $\Omega =\mathcal{C}(\mathbb{R}_+,\mathbb{R}_+)$, let $Z$ be the canonical map and $\mathcal{F}_t=\sigma\{Z_s,0\le s\le t\}$ be the canonical filtration. We denote by $\mathbb{P}^{\delta,y}$ the law of the Bessel process of dimension $\delta$ starting from $y$, $y>0$. \\
Let us state the following result from Jeanblanc, Yor and Chesney {\cite{jeanblanc-yor-chesney} (Proposition 6.1.5.1, page 364).
\begin{proposition}
The following absolute continuity relation between a 
Bessel process of dimension $\delta$ and a Bessel process of dimension $2$ holds
\[
\mathbb{P}^{\delta,y}\left|_{\mathcal{F}_t}\right.=\left( \frac{Z_t}{y}\right)^{\delta/2-1}\exp\left\{- \frac{(\delta/2-1)^2}{2}\int_0^t\frac{ds}{Z_s^2} \right\} \mathbb{P}^{2,y}\left|_{\mathcal{F}_t}\right. , \forall t\geq 0.
\]
\end{proposition}
It seems difficult to use such a relation in order to simulate Bessel hitting times for any dimension $\delta$ even if it leads to study only the $2$-dimensional case. The use of Radon-Nikodym's derivative happens to be unuseful for numerical purposes.
\section{Additivity of Bessel processes}
\par\noindent\null $\,$ \\ 

\label{sec:add}
An important property, due to Shiga and Watanabe \cite{Shiga-Watanabe}, is the additivity property for the family of squared Bessel processes usually denoted  by BESQ. Let us denote by
$\mathbb{P} _1*\mathbb{P}_2$ the convolution of $\mathbb{P}_1$ and $\mathbb{P}_2$, where  $\mathbb{P}_1$ and ${\mathbb{P}_2}$ are probability measures. In the following, we denote by $\mathbb{Q}^{\delta,x_0}$ the law of the squared Bessel process of dimension $\delta$ starting from $x_0$.
\begin{proposition} For every $x_0,\, x'_0 \geq 0$ and for every $\delta,\, \delta ' \geq 0$ we have 
\begin{equation}\label{eq:convol}
\mathbb{Q}^{\delta,x_0}*\mathbb{Q}^{\delta',x'_0}=\mathbb{Q}^{\delta+\delta',x_0+x'_0}.
\end{equation}
\end{proposition}
\section{Notations and preliminary results}
\par\noindent\null $\,$ \\ 

We start by recalling results and notations introduced in \cite{Deaconu-Herrmann} that will be needed in the follow-up.

Consider the first hitting time of a curved boundary for the Bessel process of dimension $\delta$ starting from the origin. Let $\psi(t)$ denote the boundary, and introduce the following hitting time:
\begin{equation}
\label{taupsi}
 \tau_\psi=\inf\{ t>  0;\ X^{\delta,0}_t\geq \psi(t)\}.
\end{equation}
For some suitable choice of the boundary, the distribution of $\tau_\psi$ can be explicitly computed. The result is based on the method of images (see for instance \cite{daniels_1969} for the origin of this method  and \cite{lerche_1986} for a complete presentation).
\begin{proposition}\label{prop:rappel}
Set $a>0$ and $\delta >0$. Let us consider the following curved boundary:
\begin{equation}
\label{eq:psi}
 \psi_{a,\delta}(t)=\sqrt{2t\log\ds\frac{a}{\Gamma (\frac{\delta}{2}) t^{\frac{\delta}{2}}2^{\frac{\delta}{2}-1}}},
\quad \mbox{ for }\quad t\in {\rm {Supp}} (\psi):=\left[ 0 ,\,T_{a,\delta} \right],
\end{equation}
where $T_{a,\delta}$ is the largest definition time and is given by
\begin{equation}
\label{def_time}
T_{a,\delta} := \left(\frac{a}{\Gamma(\frac{\delta}{2})2^{\frac{\delta}{2}-1}}\right)^{\frac{2}{\delta}}.
\end{equation}
We can express explicitly the distribution of $\tau_\psi$ (simplified notation corresponding to $ \tau_{\psi_{a,\delta}}$). It has its support in ${\rm {Supp}} (\psi)$ and is given by
\begin{equation}
\label{densite_tau}
\begin{array}{ll}
\mathbb{P}_0(\tau_\psi \in \dint t) &= \ds\frac{1}{2at}
\left( 
2t \log \ds\frac{a}{\Gamma (\frac{\delta}{2})t^{\frac{\delta}{2}}2^{\frac{\delta}{2}-1}}
\right)^{\frac{\delta}{2}}\ind{{\rm {Supp}} (\psi)}(t)\dint t\\
& = \ds\frac{1}{2at} (\psi_{a,\delta})^\delta (t) \ind{{\rm {Supp}} (\psi)}(t)\dint t.
\end{array}
\end{equation}
\end{proposition}
Let us note that the maximum of the function $\psi_{a,\delta}(t)$ is reached for $t_{\max}(a)= \frac{T_{a,\delta}}{e}$ and equals 
\begin{equation}
 \label{eq:maxcalcul}
\begin{array}{ll}
W_{a,\delta}& :=\ds\sup_{t\in  {\rm {Supp}} (\psi)}\ \psi_{a,\delta}(t)=\sqrt{\frac{\delta}{e}\left(\frac{a}{\Gamma(\frac{\delta}{2})2^{\frac{\delta}{2}-1}}\right)^{\frac{2}{\delta}}}= \sqrt { \frac{\delta}{e}T_{a,\delta}}.\\
\end{array}
\end{equation}
Moreover the distribution $u_{a,\delta } (t,x)\dint x:=\mathbb{P}(X_t^{\delta,0}\in\dint x,\, \tau_{\psi}>t)$ has the form
\begin{equation}
\label{density_Bessel}
u_{a,\delta} (t,x) = \left( \ds\frac{1}{2^{\frac{\delta}{2}-1} \Gamma(\frac{\delta}{2}) t^{\frac{\delta}{2}} }\exp{\left(-\ds\frac{x^2}{2t}\right)} -\ds\frac{1}{a}\right)x^{\delta-1}.
\end{equation}
{\bf Scaling property} Notice a scaling property which will be used in the sequel. Using relations \eqref{eq:psi} and \eqref{def_time} we obtain :
\begin{align}
\label{scaling}
\psi_{a,\delta}^2(t)&  = 2t\log \left(\frac{T_{a,\delta}}{t}\right)^{\frac{\delta}{2}} =\delta T_{a,\delta}\left( \frac{t}{T_{a,\delta}}\log \left(\ds\frac{T_{a,\delta}}{t}\right) \right)=\delta T_{a,\delta}\, \Phi^2 \left(\ds\frac{t}{T_{a,\delta}} \right),
\end{align}
where 
\begin{equation}
\label{Phi}
\Phi(t) = \sqrt{ t \log \left( \frac{1}{t}\right) } \ind{[0,1]}(t).
\end{equation} 
\emph{\sc Proof for Proposition \ref{prop:rappel}.}\\ This result was already presented in the particular case where the dimension of the Bessel process $\delta$ is strictly larger than $1$. In that situation, the Bessel process satisfies a SDE:
\begin{equation}\label{eq:sde}
X^{\delta,x_0}_t=x_0+\frac{\delta-1}{2}\int_0^t(X^{\delta,x_0}_s)^{-1}\dint s+B_t.
\end{equation}
The result is then linked to properties of the associated partial differential operator. Let us note that the Bessel process is always non negative.  But the Bessel processes associated with dimensions smaller than $1$ are  not semi-martingales and so, the dynamic of the processes are related to the local time at $0$ and do not satisfy \eqref{eq:sde}. For this reason we will propose a proof based only on stochastic tools (which is inspired by arguments developed for the Brownian motion by Lerche \cite{lerche_1986}). We introduce first the squared Bessel process $Y_t=(X^{\delta,x_0}_t)^2$ which satisfies the SDE for all $\delta>0$:
\begin{equation}\label{eq:carre}
Y_t=x_0^2+\delta t+2\int_0^t\sqrt{Y_s}\dint B_s.
\end{equation}
Let us denote the transition probabilities by $p_{y_0}(t,\dint x):=\mathbb{P}((X^{\delta,x_0}_t)^2\in\dint x)=\mathbb{P}_{y_0}(Y_t\in\dint x)$ with $y_0=x_0^2$. The density is given by:
 \begin{equation}
 \label{eq:densit-expr}
 p_y(t,x)=\frac{1}{2t}\left( \frac{x}{y} \right)^{\nu/2}\exp\left(-\frac{x+y}{2t}\right)I_\nu\Big(\frac{\sqrt{xy}}{t}\Big), \quad\mbox{for}\quad t>0,\ y>0,\ x\ge 0,
 \end{equation}
where $\nu=\delta/2-1$ and $I_\nu(z)$ is the Bessel function whose expression is written:
\[
I_\nu(z)=\sum_{n=0}^\infty \left(\frac{z}{2}\right)^{\nu+2n}\frac{1}{n!\Gamma(\nu+n+1)}.
\]
Moreover, for $y=0$, we get
\[
p_0(t,x)=\frac{x^{\frac{\delta}{2}-1}}{(2t)^{\delta/2}\Gamma(\delta/2)}\, e^{-\frac{x}{2t}},\quad x\ge 0.
\]
{\bf Step 1.} \\ 
Let us denote by $\mathbb{P}^{t,y}_x$ the distribution of the squared Bessel bridge starting at $x$ and reaching $y$ at time $t$ and $\mathbb{P}_x$ stands for the distribution of the squared Bessel process starting at $x$, we denote lastly $(\mathcal{F}_t)_{t\ge 0}$ the filtration associated to the Bessel process. \\
Let $t_0>0$. 
Simple computations (using the transition probabilities and the Markov property of the Bessel process) permit obtaining the following Radon-Nikodym derivative for $t<t_0$ (see, for instance, \cite{Revuz-Yor-99} page 433):
\begin{equation}\label{eq:R-N}
\left.\frac{\dint\mathbb{P}^{t_0,y}_x}{\dint\mathbb{P}_x}\right|_{\mathcal{F}_t}=\frac{p_{Y_t}(t_0-t,y)}{p_x(t_0,y)},\quad y>0,\ x\ge 0.
\end{equation}
Let us note that this formula cannot be extended to the case $t=t_0$ since $\mathbb{P}^{t_0,y}_x$ is not absolutely continuous with respect to $\mathbb{P}_x$ on $\mathcal{F}_{t_0}$. For $y=0$, the result can be obtained by continuity:
\begin{equation}\label{eq:R-N-1}
\left.\frac{\dint\mathbb{P}^{t_0,0}_x}{\dint\mathbb{P}_x}\right|_{\mathcal{F}_t}=\lim_{y\to 0}\frac{p_{Y_t}(t_0-t,y)}{p_x(t_0,y)},\ x\ge 0.
\end{equation}
Let us compute explicitly the r.h.s of \eqref{eq:R-N-1}. By \eqref{eq:densit-expr} and for $x>0$, we get
\begin{align*}
\left.\frac{\dint\mathbb{P}^{t_0,0}_x}{\dint\mathbb{P}_x}\right|_{\mathcal{F}_t}&=\lim_{y\to 0}\frac{\frac{1}{2(t_0-t)}\Big(\frac{y}{Y_t} \Big)^{\nu/2}\exp\Big(-\frac{Y_t+y}{2(t_0-t)}\Big)I_\nu\Big( \frac{\sqrt{Y_t y}}{t_0-t} \Big)}{\frac{1}{2t_0}\Big(\frac{y}{x} \Big)^{\nu/2}\exp\Big(-\frac{x+y}{2t_0}\Big)I_\nu\Big( \frac{\sqrt{x y}}{t_0} \Big)}\\
&=\left(\frac{t_0}{t_0-t}\right)^{\nu+1}\exp\left( -\frac{Y_t}{2(t_0-t)}+\frac{x}{2t_0} \right).
\end{align*}
This result can be expressed with respect to both the transition probability and the invariant measure $\mu$ satisfying $\mu(x)p_x(t,y)=\mu(y)p_y(t,x)$ that is $\mu(x)=x^\nu$. Defining
\begin{equation}
\label{eq:def:xi}
\xi(t,x)=\frac{p_0(t,x)}{\mu(x)}=\frac{1}{(2t)^{\delta/2}\Gamma(\delta/2)}\, e^{-\frac{x}{2t}},\quad\mbox{for}\ x>0, \, t>0,
\end{equation}
we obtain finally for any $t<t_0$:
\begin{equation}
\label{eq:express-RN}
D_t:=\left.\frac{\dint\mathbb{P}^{t_0,0}_x}{\dint\mathbb{P}_x}\right|_{\mathcal{F}_t}=\frac{\xi(t_0-t,Y_t)}{\xi(t_0,x)}.
\end{equation}
Let us just note that this result can be extended continuously to the case $x=0$ by defining $\xi(t,0)=\lim_{x\to 0}\xi(t,x)=(2t)^{-\delta/2}\Gamma(\delta/2)^{-1}$
and that $(D_s)_{s\le t}$ is a martingale with respect to $\mathbb{P}_x$ for $t<t_0$.\\ \\
{\bf Step 2.}\\ We prove now that $U(t_0,x)$ defined by
\[
U(t_0,x)\dint x:=\mathbb{P}_0(Y_{t_0}\in\dint x,\,\tau_{\psi^2}(Y)>t_0)
\]
satisfies
\begin{align*}
U(t_0,x)=\left(\frac{1}{(2t_0)^{\delta/2}\Gamma(\delta/2)}\, e^{-\frac{x}{2t_0}}-\frac{1}{2a}\right)x^{\nu},
\end{align*}
which directly implies \eqref{density_Bessel}.\\ 
By conditioning, we obtain
\begin{align*}
U(t_0,x)=\mathbb{P}_0(\tau_{\psi^2}(Y)>t_0|Y_{t_0}= x)p_0(t_0,x).
\end{align*}
Due to a time inversion transformation and by using the Radon-Nikodym derivative given in Step~1, the following equation yields for $x<\psi^2(t_0)$:
\begin{align*}
U(t_0,x)&=\lim_{t\to t_0,\, t<t_0}\mathbb{P}_x^{t_0,0}(\hat{\tau}(Y)>t)p_0(t_0,x)\\
&=\lim_{t\to t_0,\, t<t_0}\mathbb{E}_x[\ind{\{ \hat {\tau}>t \}}D_t]p_0(t_0,x)\\
&=\lim_{t\to t_0,\, t<t_0}\mathbb{E}_{x}[D_{t\wedge \hat{\tau}}]p_0(t_0,x)-\lim_{t\to t_0,\, t<t_0}\mathbb{E}_{x}[\ind{\{ \hat{\tau}\le t \}}D_{ \hat{\tau}}]p_0(t_0,x),
\end{align*}
where $\hat{\tau}=\tau_{\psi^2(t_0-\cdot)}$. Since $D_t$ is a continuous martingale, the optimal stopping theorem (in the time inverse filtration) leads to $\mathbb{E}_{x}[D_{t\wedge \hat{\tau}}]=D_0=1$. Moreover the function $\psi^2(t)$ has the following property: if $0\le x<\psi^2(t)$ then $\xi(t,x)<\frac{1}{2a}$, if $x>\psi^2(t)$ then $\xi(t,x)>\frac{1}{2a}$ and $\xi(t,\psi^2(t))=\frac{1}{2a}$. In other words, the stopping time $\hat{\tau}$ can be defined as follows:
\[
\hat{\tau}=\inf\Big\{t>0:\ \xi(t_0-t,Y_t)\ge \frac{1}{2a}\Big\}=\inf\Big\{t>0:\ D_t\ge \frac{1}{2a\xi(t_0,x)}\Big\}.
\]
We deduce that 
\[
\mathbb{E}_{x}[\ind{\{ \hat{\tau}\le t \}}D_{ \hat{\tau}}]=\frac{1}{2a\xi(t_0,x)}\,\mathbb{P}_x(\hat{\tau}\le t).
\]
Therefore
\begin{align*}
U(t_0,x)=p_0(t_0,x)\Big( 1-\frac{1}{2a\xi(t_0,x)}\lim_{t\to t_0}\mathbb{P}_x(\hat{\tau}\le t) \Big)=p_0(t_0,x)\Big( 1-\frac{\mathbb{P}_x(\hat{\tau}\le t_0)}{2a\xi(t_0,x)} \Big).
\end{align*}
Since $\psi$ is a continuous function with $\psi(0)=0$, we deduce that $\mathbb{P}_x(\hat{\tau}\le t_0)=1$. Thus, we obtain the result
\[
U(t_0,x)=p_0(t_0,x)\Big( 1-\frac{1}{2a\xi(t_0,x)} \Big).
\]
{\bf Step 3.}\\ As an immediate consequence, the expression of $u_{a,\delta}(t,x)$ defined by \eqref{density_Bessel} leads to 
\begin{align*}
\mathbb{P}_0(\tau_{\psi^2}(Y)>t)&=\mathbb{P}_0(\tau_\psi(X)>t)=
\int_0^{\psi^2_{a,\delta}(t)}U(t,x)\dint x=\int_0^{\psi_{a,\delta}(t)}2wU(t,w^2)\dint w\\
&=\int_0^{\psi_{a,\delta}(t)}u_{a,\delta}(t,w)\dint w.
\end{align*}
The density of the hitting time can be easily obtained by derivation, see Proposition 2.2 in \cite{Deaconu-Herrmann}.\hfill{$\Box$}

\rm
\section{Algorithm for approaching the hitting time}
\label{sec:algo}
\par\noindent\null $\,$ \\ 

In a previous paper \cite{Deaconu-Herrmann}, the authors developed an algorithm in order to simulate, in just a few steps, the Bessel hitting time. This was done for integer dimensions $\delta \geq 1$. The particular connexion between the Bessel process and the $\delta$-dimensional Brownian motion gives in this case a geometrical interpretation in terms of the exit problem from a disk for a Brownian motion in dimension $\delta$. This geometrical approach doesn't work any longer for non integer dimensions. In order to handle this difficulty, we construct here a new algorithm which is based on the additivity property expressed in Section \ref{sec:add}.\\ \\
\textbf{Some notations:} Let us start by making some notations. For a given dimension $\delta > 0$ denote $\delta ' := \delta - \lfloor \delta\rfloor$. Consider also $\gamma \in [0,1)$ (close to 1) and $L>0$. We define, for $x > 0$ 
\begin{equation}
\label{A-partie-entiere-delta}
I(\delta, x): = 2^{\lfloor \delta\rfloor/2-1}\Gamma\left( \frac{\lfloor \delta\rfloor}{2} \right)\left( \frac{\sqrt{e}\gamma(L^2-x^2)}{\sqrt{(\lfloor\delta\rfloor-\delta\gamma)x^2+
\delta\gamma L^2}+\sqrt{\lfloor\delta\rfloor x^2}} \right)^{\lfloor \delta\rfloor}
\end{equation}
and
 \begin{equation}
\label{Gammma-partie-non-entiere-delta}
N (\delta, x): = 2^{\delta'/2-1}\Gamma\left( \frac{ \delta '}{2} \right)\left( \frac{\sqrt{e}\gamma(L^2-x^2)}{\sqrt{(\delta '-\delta\gamma)x^2+
\delta\gamma L^2}+\sqrt{\delta ' x^2}} \right)^{ \delta '}.
\end{equation}
\\
\centerline{\line(1,0){490}}\\[5pt]
\emph{\textbf{Algorithm (NI) :} Simulation of $\tau_L=\inf\{t\ge 0:\ X_t^{\delta,x_0} =L\}$\\
\textbf{Initialization:} $\Theta_0=0$, $M(0)=0$, $\gamma\in (0,1)$ some parameter chosen close to $1$.\\
\mathversion{bold}
\textbf{Step $n$, ($n\ge 1$):}
\mathversion{normal}
The Bessel process starts at time $\Theta_{n-1}$ in $M(n-1)$. 
While $L^2-M^2(n-1) >\varepsilon$ do:
\begin{itemize}
\item [\mathversion{bold}
\textbf{($n.1$)}\mathversion{normal}] Construct a 
Bessel process of dimension $\lfloor \delta \rfloor$ starting from  $M(n-1)$ and stop this process at time $\theta_{n}^{(1)}$, the exit time of the $\lfloor \delta\rfloor$-dimensional Brownian motion $B_t$ from the moving sphere centered in $\underline{M}(n-1):=(M(n-1),0,0,\ldots,0)\in\mathbb{R}^{\lfloor \delta\rfloor}$ and with radius $\psi_{\alpha_n,\lfloor \delta\rfloor}(t)$, where 
\[
\label{def:An}
\alpha_n = I (\delta , M(n-1)),
\]
following the definition given in \eqref{A-partie-entiere-delta}.
\item [\mathversion{bold}
\textbf{($n.2$)}\mathversion{normal}] Construct also a second Bessel process, independent with respect to the previous one, of dimension $\delta':=\delta-\lfloor \delta\rfloor$ starting from $0$. Stop this process the first time $\theta_n^{(2)}$ it hits the curved boundary $\psi_{\beta_n,\delta'}(t)$, where
\[
\label{def:Gamma-n}
\beta_n = N (\delta , M(n-1)),
\]
following the definition given in $(\ref{Gammma-partie-non-entiere-delta})$.
\item [\mathversion{bold}
\textbf{($n.3$)}\mathversion{normal}] Define the stopping time (comparison of the two hitting times) 
\begin{equation}
\label{stopping-time}
\theta_n=\inf\{ \theta_n^{(1)}, \theta_n^{(2)}\}.
\end{equation}
First notice that the additivity property of the Bessel processes ensures that $(X_t^{\delta,M(n-1)})^2$ has the same distribution as the sum of two independent processes defined in steps \textbf{($n.1$)}\mathversion{normal} and 
\textbf{($n.2$)}\mathversion{normal}. We denote by
\[
(Z_t^{\delta,M(n-1)})^2=\Vert \underline{M}(n-1)+B_t\Vert^2+(X^{\delta',0}_t)^2.
\]
The values of $\alpha_n$ and $\beta_n$ have been chosen in order to ensure the following bound:
\begin{equation}\label{eq:majora}
\sup_{t\le \theta_n}(Z_t^{\delta,M(n-1)})^2\le
\sup_{t\le \theta_n^{(1)}}\Vert \underline{M}(n-1)+B_t\Vert^2+\sup_{t\le \theta_n^{(2)}}(X^{\delta',0}_t)^2\le M^2(n-1)+\gamma(L^2-M^2(n-1)).
\end{equation}
In particular, since $\gamma<1$,
\(
\sup_{t\le \theta_n}Z_t^{\delta,M(n-1)}< L.
\)
\\
We lastly define $M(n)=Z_{\theta_n}^{\delta,M(n-1)}< L$ and $\Theta_n=\Theta_{n-1}+\theta_n$. This achieves the $n$-th step.
\end{itemize}
 {\bf Outcome:} $N^\varepsilon$ is then the number of steps entirely completed that is the first one in the algorithm such that $L^2-(M(n))^2\le \varepsilon$, $\Theta_{N^\varepsilon}$ the approximate hitting time and $M(N^\varepsilon)$ the approximate exit position.\\
}
\centerline{\line(1,0){490}}
\noindent
Figure 1 presents several paths of the random walk $(M(n),\,n\ge 0)$ defined by the algorithm (NI) for $\delta=2.7$, $\gamma=0.9$ and for the level $L=5$.
\begin{figure}
\centerline{\includegraphics[scale=0.7]{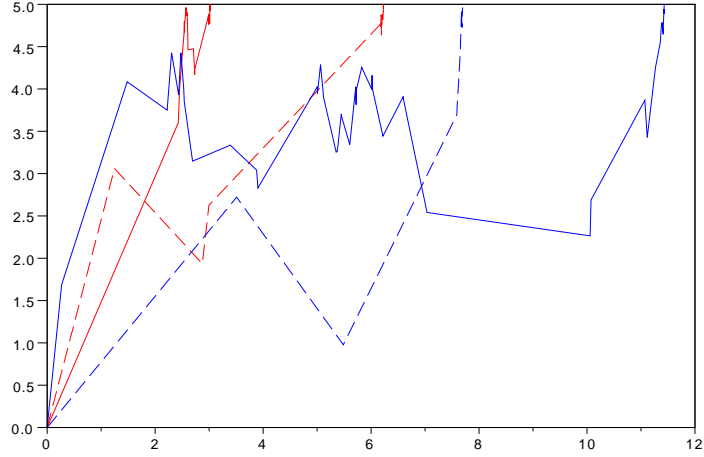}}
\caption{Several paths of the random walk}
\end{figure}
\noindent \begin{rem} 1. The upper-bound \eqref{eq:majora} will be discussed in the proof of Theorem \ref{thm:algo} (see the bound \eqref{eq:introgamm} in Step~2).\\
2. The first step of the algorithm could be simplified since the starting position is $0$. The procedure then consists in first choosing $\alpha_1$ such that $\psi_{\alpha_1,\delta}(t)<L$ for all $t\ge 0$. Then we simulate the first hitting time of this moving boundary which is linked to the Gamma distribution. We denote this random variable by $\theta_1$ and we compute the value of the process at this time: 
\[
M(1)=\psi_{\alpha_1,\delta}(\theta_1).
\]
Finally we set $\Theta_1=\Theta_0+\theta_1$. Even if this modified version of the algorithm (NI) seems simpler (just the first step is different), in the following we will prove only the convergence results associated to the algorithm (NI). \\
\end{rem} 
\subsection*{\bf Realization of the algorithm}
\par\noindent $\,$

One particular important task in this procedure is the simulation of $ Z_{\theta_n}^{\delta,M(n-1)}$ in the $n$-th step. The method we use is the following:
\begin{itemize}
\item If $\theta_n=\theta_n^{(1)}$ then 
\[
\Vert \underline{M}(n-1)+B_{\theta_n}\Vert^2= M^2(n-1)+2M(n-1)\pi_1(U)\psi_{\alpha_n,\lfloor \delta \rfloor}(\theta_n)+\psi^2_{\alpha_n,\lfloor \delta \rfloor}(\theta_n),
\]
where $\pi_1$ is the projection on the first coordinate and $U$ is a random variable in $\mathbb{R}^{\lfloor \delta\rfloor}$ uniformly distributed on the sphere of radius $1$. It suffices now to simulate $X^{\delta',0}_{\theta_n}$. Since $\theta_n^{(1)}$ and $\theta_n^{(2)}$ are independent, we get
\[
\mathbb{E}\Big[f(X^{\delta',0}_{\theta_n})\Big\vert \theta_n^{(2)}>\theta_n^{(1)}\Big]=\int_{\mathbb{R}_+}f(x)w(\theta_n,x)\,\dint x,
\]
where 
\[
w(t,x)\dint x=\mathbb{P}(X^{\delta',0}_t\in \dint x\vert \tau_\psi >t)=\frac{\mathbb{P}(X^{\delta',0}_t\in \dint x,\,\tau_\psi >t)}{\mathbb{P}( \tau_\psi >t)}=\frac{u(t,x)}{\int_0^{\psi(t)}u(t,x)\,\dint x}
\]
and $u(t,x)$ stands here for $u_{\beta_n, \delta '}(t,x)$ which was already defined in the previous paper \cite{Deaconu-Herrmann} and was restated in \eqref{density_Bessel}. More precisely
\[
u(t,x)=\left( \frac{1}{2^{\delta'/2-1}}\frac{1}{t^{\delta'/2}}\frac{1}{\Gamma(\delta'/2)}\exp\left(-\frac{x^2}{2t}\right)-\frac{1}{\beta_n} \right) x^{\delta'-1}.
\]
 Here for notational simplicity the index $\psi$ of $\tau_\psi$ stands for $\psi_{\beta_n,\delta'}$.
Let us just note that the support of $w(\theta_n,\cdot)$ is $[0,\psi_{\beta_n,\delta'}(\theta_n)]$. In order to simulate $X^{\delta',0}_{\theta_n}$, given $\tau_\psi >\theta_n$, we employ a rejection sampling method. Let $S$ be a random variable defined on the interval $[0,\psi_{\beta_n,\delta'}(\theta_n)]$ with the probability density function:
\[
r(x)=\frac{\delta'x^{\delta'-1}}{(\psi_{\beta_n,\delta'}(\theta_n))^{\delta'}},\quad \mbox{for}\quad 0\le x\le \psi_{\beta_n,\delta'}(\theta_n).
\]
This variable can be easily sampled by using a standard uniform random variable $V$: $S$ has the same distribution as
\[
\psi_{\beta_n,\delta'}(\theta_n)V^{1/\delta'}.
\]
Considering the following constant:
\[
C=\frac{(\psi_{\beta_n,\delta'}(\theta_n))^{\delta'}}{\delta'}\left( \frac{1}{2^{\delta'/2-1}\theta_n^{\delta'/2}\Gamma(\delta'/2)}-\frac{1}{\beta_n} \right), 
\]
we observe that $u(\theta_n,x)\le C r(x)$ for all $x$. Then the procedure is the following
\begin{enumerate}
\item Sample two independent r.v. $U_*$ and $S$ on respectively $[0,1]$ and $[0,\psi_{\beta_n,\delta'}(\theta_n)]$. The first one is uniformly distributed and the p.d.f. of the second one is given by $r(x)$. 
\item If $\displaystyle U_*\le \frac{u(\theta_n,S)}{Cr(S)}$  define $\xi '=S$ otherwise return to the first step.
\end{enumerate}
With this algorithm, the p.d.f. of $\xi '$ is equal to $w(\theta_n,x)$, it has the same distribution as $X^{\delta',0}_{\theta_n}$ given $\theta_n=\theta_n^{(1)}$.\\
Finally we obtain 
\[
Z_{\theta_n}^{\delta,M(n-1)}=\sqrt{(\xi ')^2+M^2(n-1)+2M(n-1)\pi_1(U)\psi_{\alpha_n,\lfloor \delta \rfloor}(\theta_n)+\psi^2_{\alpha_n,\lfloor \delta \rfloor}(\theta_n)}.
\]
\item If $\theta_n=\theta_n^{(2)}$ the result is quite similar. We obtain
\[
Z_{\theta_n}^{\delta,M(n-1)}=\sqrt{\psi^2_{\beta_n,\delta'}(\theta_n)+M^2(n-1)+2M(n-1)\pi_1(U)\xi+\xi^2 }
\]
where $\xi $ is obtained in a similar way as $\xi '$. We just replace $\delta'$ by $\delta$ and $\beta_n$ by $\alpha_n$.
\end{itemize}
\begin{thm}
\label{thm:algo}
Set $\delta  \geq 1$. The number of steps $N^\varepsilon$ of the algorithm  (NI) is almost surely finite. Moreover, there exist constants $C_\delta>0$ and $\varepsilon_0(\delta)>0$, such that
\[
\mathbb{E}[N^\varepsilon]\le C_\delta|\log\varepsilon|, \mbox{ for all }\,\, \varepsilon\le\varepsilon_0(\delta).
\]
Furthermore $\Theta_{N^\varepsilon}$ converges in distribution towards $\tau_L$, the hitting time of the level $L$ for the $\delta$ dimensional Bessel process as $\varepsilon\to 0$.
\end{thm}
The following histograms present the distribution of the hitting times $\Theta_{N^\epsilon}$ for $\delta=1.5$ and $\delta=7.5$. Of course, the higher  the dimension of the Bessel process is, the smaller  the hitting time is.\\
\centerline{\includegraphics[scale=0.5]{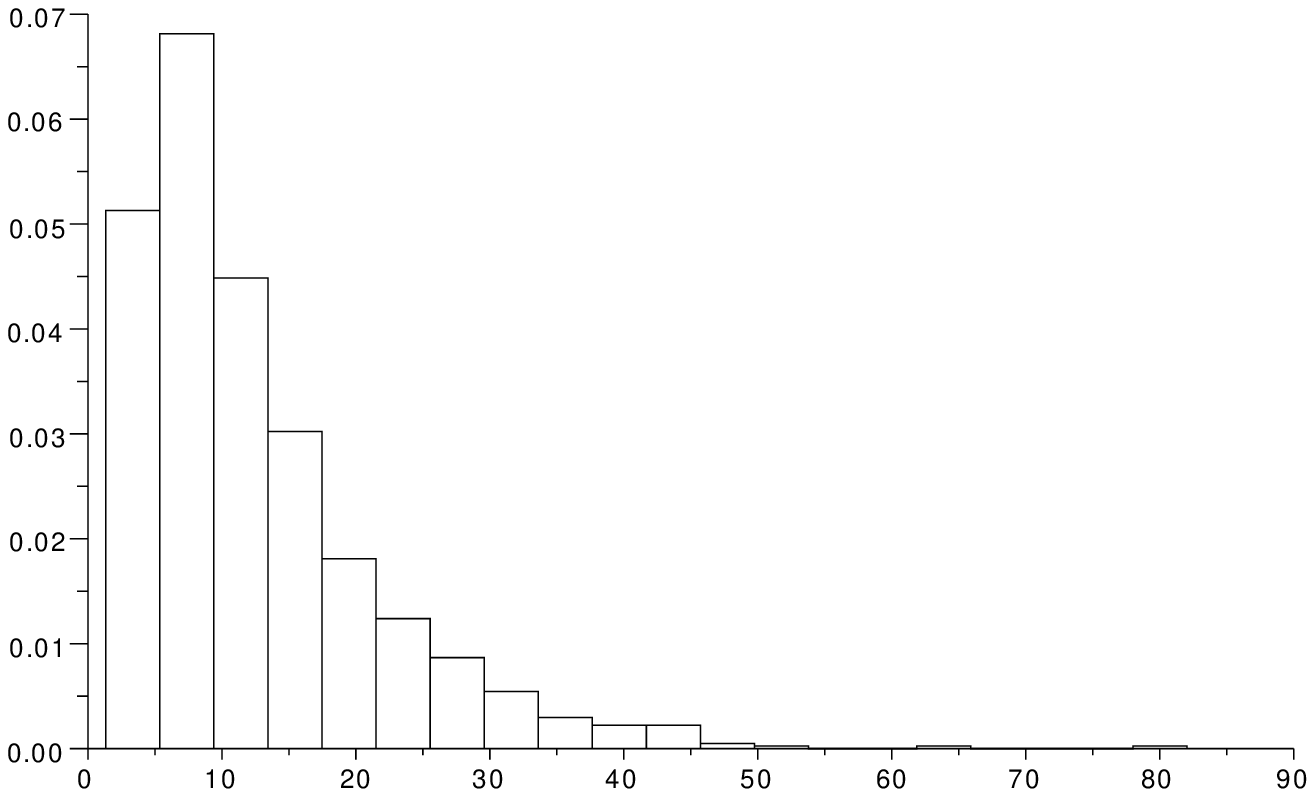}
\includegraphics[scale=0.5]{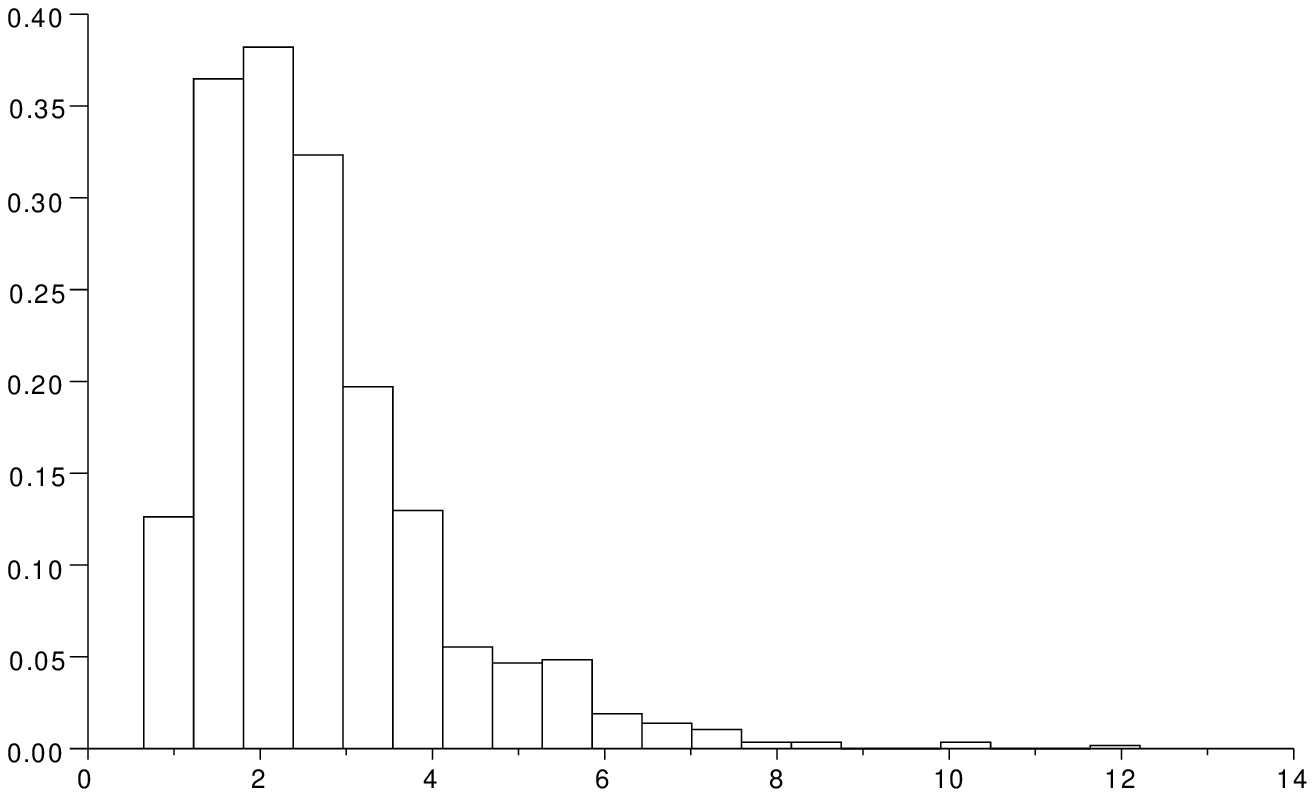}}
These simulations ($1000$ for each histogram) have been realized with the following choice of parameters: $\gamma=0.9$, $L=5$ and $\varepsilon=0.001$, $\delta =1.5$ on the left and $\delta =7.5$ on the right.
\begin{proof} Instead of considering the Markov chain $(M(n),\, n\ge 0)$, we focus our attention to the squared process $R(n)=(M(n))^2$ and we stop the algorithm as soon as $R(n)$ becomes larger than $L^2-\varepsilon$.\\ \\
\emph{\bf Step 1. Definition and decomposition of the operator $P$} \\
We  estimate first the number of steps. Since $(R(n),\, n\ge 0)$ is an homogeneous Markov chain, let us start by computing the transition probabilities associated to $R(n)$. We introduce the operator $Pf$ defined, for any non-negative measurable function $f:\mathbb{R}_+\to\mathbb{R}_+$:
\begin{equation}
\label{operator}
Pf(x):=\mathbb{E}[f(R(n))|R(n-1)=x]=\mathbb{E}\Big[f\Big( \Vert\underline{\sqrt{x}}+B_{\theta_n}\Vert^2+(X^{\delta',0}_{\theta_n})^2 \Big)\Big],
\end{equation}
where $\theta_n$ is defined in \eqref{stopping-time}.
Since $R(n)$ is an homogeneous Markov chain, the transition $Pf$ does not depend on the time $n$. For notational simplicity we neglect some indexes: the step $n$; $X^{\delta',0}$ is replaced by $X$; $\theta_n^{(i)}$ by $\theta^{(i)}$, for $i=1,2$; $\alpha_n$ is replaced by $\alpha$ and $\beta_n$ is replaced by $\beta$.
We can express (\ref{operator}) by splitting it into two parts $P_if(x)$, $i=1,2$ with
\[
P_if(x)=\mathbb{E}\Big[f\Big( \Vert\underline{\sqrt{x}}+B_{\theta^{(i)}}\Vert^2+(X_{\theta^{(i)}})^2 \Big)\ind{\{ \theta=\theta^{(i)} \}}\Big].
\]
Thus 
\begin{equation}
\label{operator:dec}
Pf(x) = P_1f(x)+P_2f(x).
\end{equation}
The class of functions that will be considered in the following satisfies the hypothesis :\\ \\
\textbf{(H)} The function $f$ is such that 
\begin{equation*}
f^{(p)} \leq 0, \, \forall p\in\{1,2,3,4\},
\end{equation*}
this means that the first four derivatives of $f$ are negative. An example of such a function will be used later on.\\ \\
\emph{\bf Step 2. Distributions associated with $\theta^{(i)}$.}\\ 
Let us denote by $U_i(t)\dint t$ the distribution of $\theta^{(i)}$ and its support $[0,s_i]$ . With the notation in \eqref{def_time} we have
\begin{equation}
\label{S1S2}
s_1= T_{\alpha,\lfloor\delta\rfloor} \mbox { and } s_2= T_{\beta, \delta '}.
\end{equation}

These distributions are those of stopping times corresponding to the function $\psi_{a, \delta}(t)$, see \eqref{eq:psi}, with a suitable parameter $a$ and a suitable dimension. By applying the scaling property \eqref{scaling} we get
\begin{equation}
\psi_{\alpha,\lfloor\delta\rfloor}(t)=\sqrt{\lfloor \delta\rfloor \, s_1 }\Phi \left( \ds\frac{t}{s_1}\right).
\end{equation}
From the same arguments, the scaling property associated to $\theta^{(2)}$ is
\[
\psi_{\beta,\delta'}(t)=\sqrt{\delta ' \, s_2}\Phi \left( \ds\frac{t}{s_2}\right), \quad\mbox{with}\ \delta'=\delta-\lfloor \delta\rfloor.
\]
Since
\begin{align*}
U_1(t)\dint t &= \mathbb{P}_0(\theta^{(1)} \in \dint t)\\
&=\frac{1}{2\Gamma(\lfloor \delta\rfloor/2)2^{\lfloor \delta\rfloor/2-1}s_1^{\lfloor\delta\rfloor/2}t}\left(2t\log\left(\frac{s_1}{t}\right)^{\lfloor \delta\rfloor /2} \right)^{\lfloor \delta\rfloor/2}\ind{[0,s_1]}(t)\dint t\\
&=\frac{c_{\lfloor\delta\rfloor}}{t}\Phi^{\lfloor \delta\rfloor}\left(\frac{t}{s_1}\right)\dint t,
\end{align*}
with $c_d=\frac{1}{\Gamma(d/2)}(\frac{d}{2})^{d/2}$. 
We therefore obtain the scaling property:
\begin{equation}
\label{eq:scal}
U_1(s_1t)=\frac{c_{\lfloor\delta\rfloor}}{s_1t}\Phi^{\lfloor \delta\rfloor}(t)\quad\mbox{and}\quad U_2(s_2t)=\frac{c_{\delta'}}{s_2t}\Phi^{\delta'}(t).
\end{equation}
With our choice of $\alpha $ and $\beta$, we have 
\begin{equation}
\label{supp:U1}
\begin{array}{ll}
s_1 & = T_{\alpha,{\lfloor \delta\rfloor}} =\left[\ds\frac{\alpha}{\Gamma(\frac{\lfloor \delta\rfloor }{2})2^{\frac{\lfloor \delta\rfloor}{2}-1}}\right]^{\frac{2}{\lfloor \delta\rfloor}}\\ 
& =\ds \frac{e\gamma^2(L^2-x)^2}{(\sqrt{(\lfloor\delta\rfloor-\delta\gamma)x+
\delta\gamma L^2}+\sqrt{\lfloor\delta\rfloor x})^2}\\
& = s_2.
\end{array}
\end{equation}
By using \eqref{eq:maxcalcul} we get
\begin{equation}
\label{supp:1}
W_{\alpha,{\lfloor \delta\rfloor}}=\sup_{t\in[0,s_1]}\ \psi_{\alpha,\lfloor\delta\rfloor}(t)= \sqrt{\frac{\lfloor \delta\rfloor}{e} s_1}
\end{equation}
and the same property holds for $s_2$
\begin{equation}
\label{supp:2}
W_{\beta,\delta '}=\sup_{t\in[0,s_2]}\ \psi_{\beta, \delta '}(t)= \sqrt{\frac{\lfloor \delta '\rfloor}{e} s_2}.
\end{equation}

Since we stop the Markov chain $R(n)$ in order to not hit $L^2$, we have
\begin{equation}
\label{eq:introgamm}
\begin{array}{ll}
(\sqrt{x}+W_{\alpha,\lfloor \delta\rfloor})^2+ W_{\beta,\delta '}^2 & = x+\ds\frac{\lfloor \delta\rfloor}{e} s_1+ 2\sqrt{x}\sqrt{\ds\frac {\lfloor \delta\rfloor}{e}} s_1 +\ds\frac{\delta '}{e}s_2\\
&= x+\ds\frac{\delta}{e}s_1 + 2\sqrt{x}\sqrt{\ds\frac {\lfloor \delta\rfloor}{e}} s_1\\
& =x+\gamma(L^2-x),\\
\end{array}
\end{equation}
by the choice of $s_1=s_2$ and the value of $s_1$ given in \eqref{supp:U1}.
We deduce that for $x$ close to $0$, $s_1$ is close to $\gamma eL^2/\delta$ and for $x$ close to $L^2$, $s_1$ is of the same order as
\(
e\gamma^2(L^2-x)^2/(4\lfloor \delta\rfloor L^2).
\)\\  \\
\emph{\bf Step 3. Computation of $P_1f(x)$.}\\
Using the definition of $P_1f$, we get 
\begin{equation}
\label{eq:devPf}
\begin{array}{ll}
P_1f(x)&=\mathbb{E}\left[ f\left(\Vert\underline{\sqrt{x}}+B_\theta\Vert^2+(X_\theta)^2\right)\ind{\{\theta=\theta^{(1)}\}}\right]\\ \\
&=\ds\int_{\mathcal{S}^1}\ds\int_0^{s_1}\mathbb{E}\left[ f\left( ( \sqrt{x}+\pi_1(z)\psi_{\alpha,\lfloor \delta\rfloor}(t))^2+\psi_{\alpha,\lfloor \delta\rfloor}^2(t)(1-\pi_1^2(z))+X_t^2\right)\ind{\{\theta^{(2)}>t\}}\right]U_1(t)\dint t \sigma(\dint z)\\ \\
&=\ds\int_{\mathcal{S}^1}\ds\int_0^{s_1}\mathbb{E}\left[ f\left( x+2\sqrt{x}\pi_1(z)\psi_{\alpha,\lfloor \delta\rfloor}(t)+\psi_{\alpha,\lfloor \delta\rfloor}^2(t)+X_t^2\right)\ind{\{\theta^{(2)}>t\}}\right]U_1(t)\dint t \sigma(\dint z).\\
\end{array}
\end{equation}
We denote here by $\mathcal{S}^1$ the unit sphere in $\mathbb{R}^{\lfloor \delta\rfloor }$, $\sigma(\dint z)$ the uniform surface measure on this sphere and $\pi_1 (z)$ the projection on the first coordinate and $\alpha$ is chosen so that \eqref{supp:U1} is satisfied. Let us just note that the variable $y=x+2\sqrt{x}\pi_1(z)\psi_{\alpha,\lfloor \delta\rfloor}(t)+\psi_{\alpha,\lfloor \delta\rfloor}^2(t)+X_t^2$ always stays in the interval $[0,L^2]$ on the event $\theta^{(2)}>t$. Consider the Taylor expansion of $f$ in a neighborhood of $x$. By using \textbf{(H)} we have that $f^{(4)}(t)\le 0$ on the whole interval $[0,L^2]$. Hence
\begin{equation}
\label{eq:devprinc}
P_1f(x)\le f(x)G_1+f'(x)G_2+\frac{1}{2}\,f''(x)G_3+\frac{1}{6}\,f^{(3)}(x)G_4.
\end{equation}
It suffices to compute $G_k$ for $k=1,2,3,4$ where $G_k$ is defined by 
\begin{equation}
\label{dev:Taylor}
G_k:=\int_{\mathcal{S}^1}\int_0^{s_1}\mathbb{E}\left[ \left( \psi_{\alpha,\lfloor \delta\rfloor}^2(t) + 2\sqrt{x}  \pi_1(z)\psi_{\alpha,\lfloor \delta\rfloor}(t)+X_t^2\right)^{k-1}\ind{\{\theta^{(2)}>t\}}\right]U_1(t)\dint t \sigma(\dint z),\ , \forall k=1,2,3,4.
\end{equation}
In particular $G_1=\mathbb{P} [\theta^{(2)}> \theta^{(1)}]$. 
Using symmetry arguments, the term associated to the projection vanishes, and we can split $G_2$ into two parts, $G_2=G_{2,1}+G_{2,2}$:
\[
G_{2,1}:=\int_0^{s_1}\psi_{\alpha,\lfloor \delta\rfloor}^2(t)\mathbb{P}( \theta^{(2)}>t)U_1(t)\dint t \quad\mbox{and}\quad G_{2,2}:=\int_0^{s_1}\mathbb{E}\left[ X_t^2\ind{\{\theta^{(2)}>t\}}\right]U_1(t)\dint t .
\]
By changing the variable $s_1 u=t$ and the scaling properties developed in Step 2, we obtain:
\begin{align*}
G_{2,1}&=\int_0^1 \psi_{\alpha,\lfloor \delta\rfloor}^2(s_1 u)\mathbb{P}(\theta^{(2)}>s_1 u)U_1(s_1 u)s_1 \dint u\\
&=\int_0^1 s_1 {\lfloor \delta\rfloor} \Phi^2(u)\left(\int_{u}^{1}U_2(s_1w)\dint w\right)\frac{c_{\lfloor \delta\rfloor}}{s_1u}\,\Phi^{\lfloor \delta\rfloor}(u)s_1 \dint u.
\end{align*}
Using $s_1=s_2$ and a change of variable, we get
\begin{align}
\label{eq:T21}
G_{2,1}=s_1 {\lfloor \delta\rfloor} c_{\lfloor \delta\rfloor}\int_0^1\int_u^1 \frac{1}{u}\Phi^{\lfloor\delta\rfloor+2}(u)s_1U_2(s_1v)\dint v\dint u=s_1\kappa_{2,1}
\end{align}
where
\begin{equation}
\label{def:kappa21}
\kappa_{2,1}:={\lfloor \delta\rfloor} c_{\lfloor \delta\rfloor}c_{\delta'}\int_0^1\int_u^1 \frac{1}{uv}\Phi^{\lfloor\delta\rfloor+2}(u)\Phi^{\delta'}(v)\dint v\dint u>0.
\end{equation}
Notice that $\kappa_{2,1}$ is a constant which only depends on the dimension $\delta$. 
We can also prove that there exists a constant $\kappa_{2,2}$ independent of $x$ such that $G_{2,2}=s_1 \kappa_{2,2}$. Indeed
\begin{align*}
G_{2,2}&= \int_0^{s_1}\mathbb{E}\left[ X_t^2 \ind{\{\theta^{(2)}>t\}}\right]U_1(t)\dint t=s_1\int_0^1\mathbb{E}\left[ X_{s_1 u}^2 \ind{\{\theta^{(2)}>s_1 u\}}\right]U_1(s_1u)\dint u\\
&=\int_0^1\mathbb{E}\left[ X_{s_1 u}^2 \ind{\{\forall r\le s_1u:\ X_r\le \psi_{\beta,\delta '}(r) \}}\right]\frac{c_{\lfloor \delta\rfloor}}{u}\,\Phi^{\lfloor \delta\rfloor}(u)\dint u.
\end{align*}
Using the scaling property of the Bessel process, we get
\[
G_{2,2}=\int_0^1\mathbb{E}\left[ s_1X_{u}^2 \ind{\{\forall r\le u:\ X_r\le \psi_{\beta,\delta'}(s_1r)/\sqrt{s_1} \}}\right]\frac{c_{\lfloor \delta\rfloor}}{u}\,\Phi^{\lfloor \delta\rfloor}(u)\dint u=s_1\kappa_{2,2}
\]
since $s_1=s_2$ and thus $\psi_{\beta,\delta'}(s_1r)/\sqrt{s_1}=\sqrt{\delta '} \Phi(r)$ does not depend on $x$ but only on $\delta$. To sum up, we have proved the existence of two constants $\kappa_{2,i}$, $i=1,2$, independent of $x$ satisfying
\begin{equation}
\label{eq:bilanT2}
G_2=\kappa_2  s_1, \, \rm { where }\,\, \kappa_2= \kappa_{2,1}+\kappa_{2,2}.
\end{equation}
Let us now focus our attention on $G_3$ defined in \eqref{dev:Taylor}.\\
While developing the square of $\psi_{\alpha,\lfloor \delta\rfloor}^2(t)+2\sqrt{x}\pi_1(z)\psi_{\alpha,\lfloor \delta\rfloor}(t)+X_t^2 $, we obtain $6$ terms:
\[
H_1=\psi_{\alpha,\lfloor \delta\rfloor}^4(t),\quad H_2=4x\pi_1^2(z)\psi_{\alpha,\lfloor \delta\rfloor}^2(t),\quad
H_3=X_t^4
\]
\[
H_4=4\psi_{\alpha,\lfloor \delta\rfloor}^3(t)\sqrt{x}\pi_1(z),\quad 
H_5=4\sqrt{x}\pi_1(z)\psi_{\alpha,\lfloor \delta\rfloor}(t)X_t^2, \quad
H_6=2\psi_{\alpha,\lfloor \delta\rfloor}^2(t)X_t^2.
\]
Therefore $G_3$ can be split into 6 terms: $G_3=\sum_{j=1}^6G_{3,j}$ with
\[
G_{3,j}:=\int_{\mathcal{S}^1}\int_0^{s_1}\mathbb{E}\Big[H_j\ind{\{ \theta^{(2)}>t \}} \Big]U_1(t)\dint t \,\sigma(\dint z).
\]
Now, let us compute $G_{3,j}$ for $j=1,\ldots,6$. First we note that, due to symmetry properties of the variable $z$, $G_{3,4}=G_{3,5}=0$. By similar arguments as those included in the computation of $G_2$, we get:
\begin{equation}\label{eq:T31}
G_{3,1}=s_1^2\kappa_{3,1}\quad\mbox{with}\quad \kappa_{3,1}:= {\lfloor \delta\rfloor}^2 c_{\lfloor \delta\rfloor}c_{\delta'}\int_0^1\int_u^1 \frac{1}{uv}\Phi^{\lfloor\delta\rfloor+4}(u)\Phi^{\delta'}(v)\dint v\dint u>0.
\end{equation}
\begin{equation}
\label{eq:T32}
G_{3,2}=xs_1\kappa_{3,2}\quad\mbox{with}\quad \kappa_{3,2}:=4\kappa_{2,1}\int_{\mathcal{S}^1}\pi_1^2(z)\sigma(\dint z)>0.
\end{equation}
\begin{equation}
\label{eq:T33}
G_{3,3}=s_1^2\kappa_{3,3}\quad\mbox{with}\quad \kappa_{3,3}:=\int_0^1\mathbb{E}\left[ X_{u}^4 \ind{\{\forall r\le u:\ X_r\le \sqrt{\delta '}\Phi(r) \}}\right]\frac{c_{\lfloor \delta\rfloor}}{u}\,\Phi^{\lfloor \delta\rfloor}(u)\dint u>0.
\end{equation}
\begin{equation}
\label{eq:T36}
G_{3,6}=s_1^2\kappa_{3,6}\quad\mbox{with}\quad \kappa_{3,6}:=2\int_0^1\mathbb{E}\left[ X_{u}^2 \ind{\{\forall r\le u:\ X_r\le  \sqrt{\delta '}\Phi(r) \}}\right]\frac{c_{\lfloor \delta\rfloor}}{u} \,{\lfloor \delta\rfloor}\, \Phi^{\lfloor \delta\rfloor+2}(u)\dint u>0.
\end{equation}
To sum up, there exist two positive constants $\kappa_3$ and $\tilde{\kappa}_3$ independent of $x$ such that
\begin{equation}
\label{eq:bilanT3}
G_3=\kappa_3 xs_1+\tilde{\kappa}_3 s_1^2.
\end{equation}
Finally, we consider the expression $G_4$ defined by \eqref{dev:Taylor}. We are not going to compute it explicitly as we have just done for the first terms (the proof would become quite boring...). 
Due to the symmetry property of the variable $z$, the expansion of $G_4$ coupled  to the computation of $(\psi_{\alpha,\lfloor \delta\rfloor}^2(t)+2\sqrt{x}\pi_1(z)\psi_{\alpha,\lfloor \delta\rfloor}(t)+X_t^2 )^3$ leads to terms which are either positive or equal to $0$. Hence
\begin{equation}
\label{eq:bilanT4}
G_4\ge 0.
\end{equation}
\emph{\bf Step 4. Computation of $P_2f(x)$.}\\
Using the definition of $P_2f$, we obtain:
\begin{align*}
P_2f(x)&=\mathbb{E}\left[ f\Big(\Vert \underline{\sqrt{x}}+B_{\theta}  \Vert^2+(X_{\theta})^2\Big)\ind{\{\theta
=\theta^{(2)}
\}} \right]\\
&=\int_0^{s_2}\mathbb{E}\left[ f\Big(\Vert \underline{\sqrt{x}}+B_t \Vert^2+\psi_{\beta,\delta'}^2(t)\Big)\ind{\{\theta^{(1)}>t\}} \right]U_2(t)\dint t\\
&=\int_0^{s_2}\mathbb{E}\left[ f\Big(x+2\langle\underline{\sqrt{x}},B_t\rangle +\Vert B_t\Vert^2+\psi_{\beta,\delta'}^2(t)\Big)\ind{\{\theta^{(1)}>t\}} \right]U_2(t)\dint t.
\end{align*}
By hypothesis \textbf{(H)},  $f''$ is non positive on the support of the Markov chain. Then, by using a Taylor expansion, the following bound holds: for any $x\in[0,L^2-\varepsilon]$
\[
P_2f(x)\le f(x)\mathbb{P}(\theta^{(1)}>\theta^{(2)})+f'(x)\int_0^{s_2}\mathbb{E}\left[\Big( 2\langle\underline{\sqrt{x}},B_t\rangle +\Vert B_t\Vert^2+\psi_{\beta,\delta'}^2(t)\Big)\ind{\{\theta^{(1)}>t\}} \right]U_2(t)\dint t.
\]
The integral expression contains three distinct terms. The first one associated to the scalar product is equal to zero since the distribution of $B_t$ given $\{\theta^{(1)}>t\}$ is rotationally invariant. The second and third terms are positive. We therefore deduce that for any function $f$ satisfying \textbf{(H)}, the following bound holds:
\begin{equation}
\label{eq:boundd}
P_2f(x)\le f(x)\mathbb{P}(\theta^{(1)}>\theta^{(2)}).
\end{equation} 
\emph{\bf Step 5. Application to a particular function $f$.}\\
Let us introduce the function $f_\varepsilon:[0,L^2)\mapsto \mathbb{R}$ defined by
\begin{equation}
\label{eq:def:func}
f_\varepsilon(x)=\log\left(\frac{L^2-x}{(1-\gamma)\varepsilon}\right),
\end{equation}
where $\gamma<1$ is the constant close to $1$ already introduced in \eqref{eq:introgamm}. Let us assume that the Markov chain $M$ starts with the initial value $x\in [0,L^2-\varepsilon]$ that is $R(0)=x$. We will prove that $f_\varepsilon(R(1))$ is a non-negative random variable. Indeed for any $y$ in the support of $R(1)$ we have:
\[
0\le y\le (\sqrt{x}+W_{\alpha,\lfloor \delta\rfloor})^2+W_{\beta,\delta '}^2,
\]
where $W_{\alpha,\lfloor \delta\rfloor}$ is defined by \eqref{supp:1} and $W_{\beta,\delta '}$ by \eqref{supp:2}. Using the identity \eqref{eq:introgamm}, we obtain
\[
f_\varepsilon(y)\ge \log\left(\frac{(L^2-x)}{\varepsilon}\right)
\ge 0
\]
since $L^2-x\ge \varepsilon$. We deduce therefore that $f_\varepsilon$ is a non-negative function on the support of the Markov chain stopped at the first exit time of the interval $[0,L^2-\varepsilon]$.

Let us now apply the operator $P$ defined in \eqref{operator} to the function $f_\varepsilon$.  Since
\[
f'_\varepsilon(x)=-\frac{1}{L^2-x}<0,\quad  f^{(k)}_\varepsilon(x)=-\frac{(k-1)!}{(L^2-x)^k}<0, \ \mbox{for}\  k=2,3,4
\]
the condition \textbf{(H)} is satisfied and we obtain by \eqref{eq:devprinc} and the computation of $G_1,\ldots, G_4$:
\begin{equation}
\label{eq:appli1}
P_1f_\varepsilon(x)\le f_\varepsilon(x)\mathbb{P}(\theta^{(2)}>\theta^{(1)})-\frac{s_1}{(L^2-x)^2}\,\Big( \kappa_2(L^2-x)+\kappa_3 x \Big),
\end{equation}
by using  \eqref{eq:bilanT2} and  \eqref{eq:bilanT3}.

Due to the explicit expression of $s_1$ in \eqref{S1S2}, there exists a constant $\kappa>0$ independent of $x$ and $\varepsilon$ such that
\begin{equation}
\label{eq:appli2}
P_1f_\varepsilon(x)\le f_\varepsilon(x)\mathbb{P}(\theta^{(2)}>\theta^{(1)})-\kappa,\quad \forall x\in[0,L^2-\varepsilon).
\end{equation}
Combined with \eqref{eq:boundd}, we get
\[
Pf_\varepsilon(x)-f_\varepsilon(x)\le -\kappa,\quad\forall x\in[0,L^2-\varepsilon].
\] 
Due to a comparison theorem of the classical potential theory, see Norris \cite{Norris-97}, we deduce that 
\begin{equation}
\label{eq:fin}
\mathbb{E}_x[N^\varepsilon]\le \frac{f_\varepsilon(x)}{\kappa}.
\end{equation} 
In particular, for $x=0$, the announced result in the statement of the theorem is proved.\\ \\
\emph{\bf Step 6. The time $ \Theta_{N^\varepsilon}$ given by the algorithm is close to the first hitting time $\tau_L$}.\\
Let us denote by $F$ (resp. $F^\varepsilon$) the cumulative distribution function of the random variable $\tau_L$ (resp. $\Theta_{N^\varepsilon}$). We construct these two random variables on the same paths ; the law of the Bessel process of dimension $\delta$ is realized as a sum of two independent Bessel processes of dimension $\lfloor \delta\rfloor$ and $\delta'$ on random time intervals until the hitting time $\Theta_{N^\varepsilon}$ and afterwards, the paths are generated just by a Bessel process of dimension $\delta$ starting in $M(N^\varepsilon)$. Since $\tau_L\ge \Theta_{N^\varepsilon}$ a.s. we immediately obtain the first bound
\[
F(t)\le F^\varepsilon(t),\quad t\ge 0.
\]  
Moreover by similar arguments as those presented in Theorem 2.9 \cite{Deaconu-Herrmann}, for $\alpha>0$, we get
\begin{align}\label{eq:upperb1}
1-F(t)&=\mathbb{P}(\tau_L>t)=\mathbb{P}(\tau_L>t,\ \Theta_{N^\varepsilon}\le t-\alpha)+\mathbb{P}(\tau_L>t,\ \Theta_{N^\varepsilon}> t-\alpha)\nonumber\\
&\le F^\varepsilon(t-\alpha)
\sup_{y\in[\sqrt{L^2-\varepsilon},L]}\mathbb{P}_y(\tau_L>\alpha)+1-F^\varepsilon(t-\alpha).
\end{align}
Applying Shiga and Watanabe's result \eqref{eq:convol}, the Bessel process of dimension $\delta>1$ is stochastically larger than the Bessel process of dimension $1$ which has the same law as $\vert B_t\vert$. Here $B$ stands for a $1$-dimensional Brownian motion. Which is why the following upper bound holds, for any $y\in [\sqrt{L^2-\varepsilon},L]$:
\begin{align}\label{eq:upperb2}
\mathbb{P}_y(\tau_L>\alpha)&\le
\mathbb{P}_{L^\varepsilon}(\tau_L>\alpha)\le\mathbb{P}_{L^\varepsilon}(\sup_{0\le t \le \alpha}|B_t|< L),\quad\mbox{with}\quad L^\varepsilon=\sqrt{L^2-\varepsilon}\nonumber\\
&\le \mathbb{P}_{L^\varepsilon}(\sup_{0\le t \le \alpha}B_t< L)\le \mathbb{P}_{0}(\sup_{0\le t \le \alpha}B_t< L-L^\varepsilon)\nonumber\\
&\le \mathbb{P}_{0}\left(\sup_{0\le t \le \alpha}B_t< \ds\frac{\varepsilon}{L}\right)\le \frac{\varepsilon}{L\sqrt{2\alpha\pi}}.
\end{align}
By combining \eqref{eq:upperb1} and \eqref{eq:upperb2}, we have
\begin{equation}\label{eq:doubleineq}
F^\varepsilon(t-\alpha)\Big( 1-\frac{\varepsilon}{L\sqrt{2\alpha\pi}} \Big)\le  F(t)\le F^\varepsilon(t),\quad t\ge 0.
\end{equation}
Consequently, $\Theta_{N^\varepsilon}$ converges to $\tau_L$ in distribution as $\varepsilon$ goes to $0$.
This ends the proof.
\end{proof}
\section{Numerical results}
\par\noindent\null $\,$ \\ 

In this section, we will discuss some numerical experiences based on the algorithm for approaching the hitting time (developed in Section \ref{sec:algo}). A particularly important task in such an iterative method is to estimate the number of steps or even the number of times the uniform random generator is used. The algorithm (NI) presented in Section \ref{sec:algo} allows to simulate hitting times for Bessel processes of non integer dimensions $\delta>1$. We will therefore only present experiences in that context and refer to the previous work \cite{Deaconu-Herrmann} for Bessel processes of integer dimensions. \\
\subsection{Number of steps versus $\varepsilon$}
\par\noindent $\,$

 The number of steps of the algorithm is of prime interest. Classical time splitting method in order to simulate particular paths of stochastic processes can be used for classical diffusion processes with regular diffusion and drift coefficients if the study is restricted to some fixed time intervals. Here the diffusion is singular, the classical methods could not be applied, nevertheless the approximation procedure (NI) developed in Section \ref{sec:algo} is of a different kind and holds at any given time. That's why we are not able to compare different methods but we will just describe the relevance of (NI) by the estimation of the average number of steps. The algorithm used in order to simulate the hitting time of the level $L$ by the Bessel process, permits obtaining an approximated hitting time $\Theta_{N^\varepsilon}$ and the corresponding position $M(N^\varepsilon)$ which satisfies:
\[
L^2-(M(N^\varepsilon))^2\le\varepsilon.
\]
The number of iterations will decrease with respect to the parameter $\varepsilon$. The average number of steps $\mathbb{E}[N^\varepsilon]$ is upper-bounded by the logarithm of $\varepsilon$ up to a multiplicative constant (Theorem \ref{thm:algo}). Let us therefore choose different values of $\varepsilon$ and approximate through a law of large number this average (we denote by $\aleph$ the number of independent simulations; here $\aleph=1000$).
\[
\varepsilon_k=0.5^k,\quad k=1,\ldots,15.
\]
 The experiences concern two different dimensions for the Bessel process $\delta=2.2$ and $\delta=4.7$ and we fix the parameter $\gamma=0.95$ appearing in the algorithm (NI) (this parameter is fixed for the whole numerical section). Note that each step of the (NI) algorithm is associated to a comparison between particular hitting times $\theta_n^{(1)}$ and $\theta_n^{(2)}$, the first one is associated with a Bessel process of dimension $\lfloor \delta\rfloor$ and the second one is associated with a Bessel process of dimension $\delta-\lfloor \delta\rfloor$. In the following we are also interested in the average number of steps $N^\varepsilon_{\rm integer}$ satisfying 
 \[
N^\varepsilon_{\rm integer}:=\#\{1\le n\le N^\varepsilon: \theta_n^{(1)}< \theta_n^{(2)}\}.
 \]
The figures represent both the estimated average number of steps $\mathbb{E}[N^\varepsilon]$ and its confidence interval (three upper curves) and the estimated average number $\mathbb{E}[N^\varepsilon_{\rm integer}]$ and its $95\%$-confidence interval (three lower curves).\\
 \centerline{\includegraphics[scale=0.5]{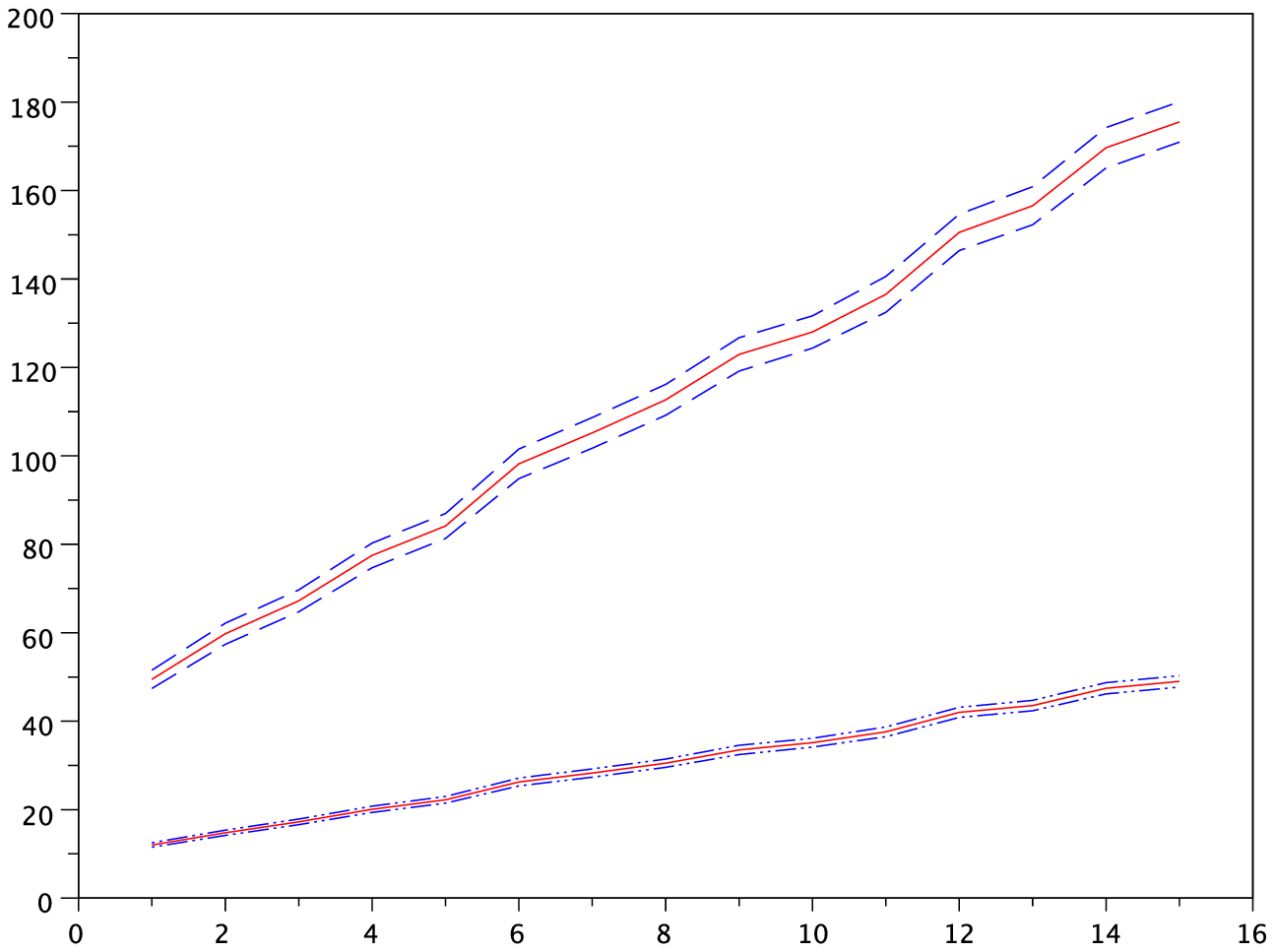}\hspace*{0.4cm}
 \includegraphics[scale=0.5]{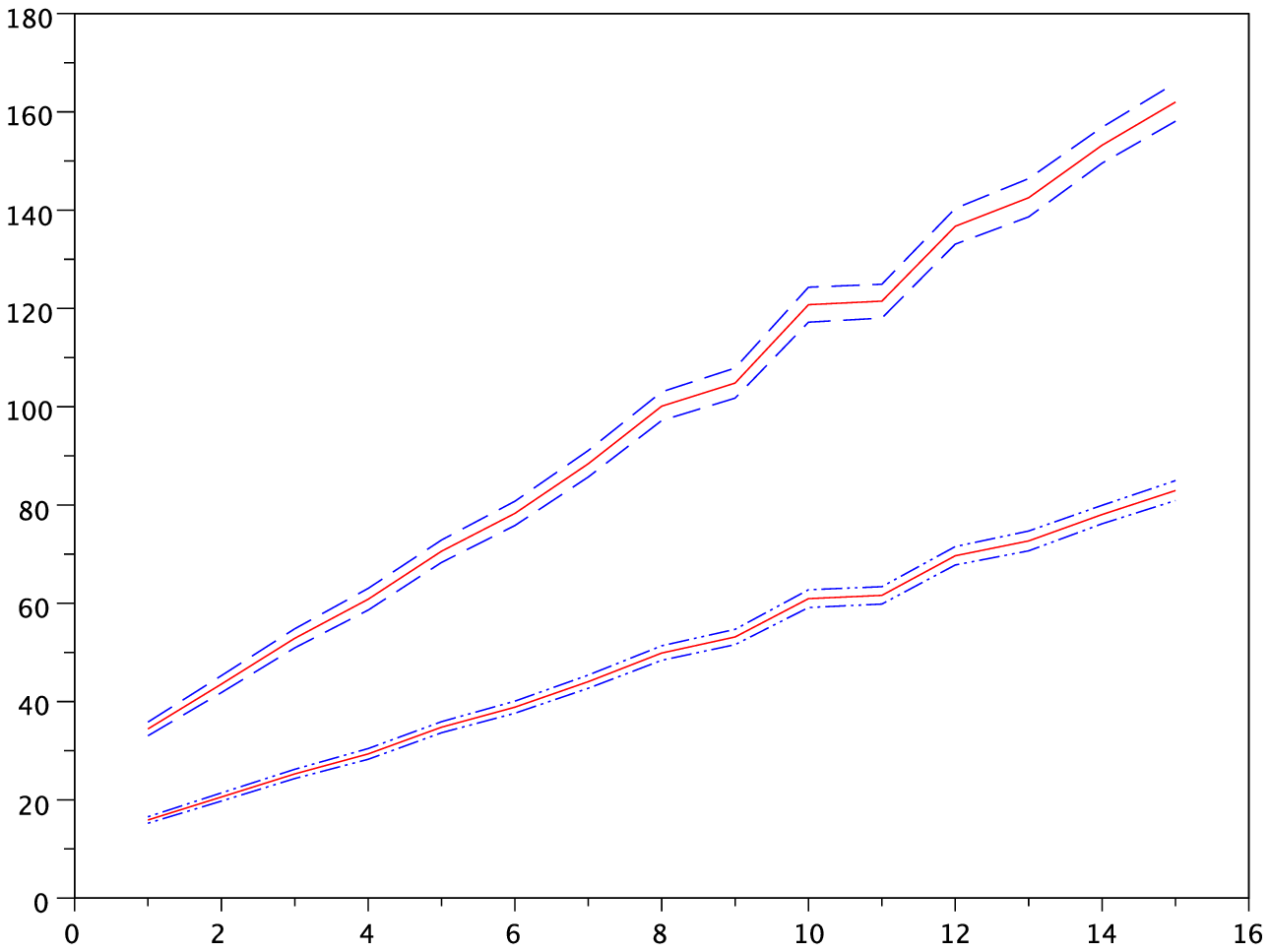}}
\centerline{Average number of steps versus $k$ for $\delta=2.2$ and $\delta=4.7$}
\subsection{Number of steps versus the dimension of the Bessel process}
\par\noindent $\,$

In \cite{Deaconu-Herrmann}, the authors pointed out that the number of steps increases as the dimension of the integer Bessel process becomes larger. Let us focus our attention on the non-integer case. We observe some surprising effects in respect to the dimension: on one hand if $\lfloor\delta\rfloor$ is fixed and the dimension increases then the average number of steps decreases, on the other hand if $\delta-\lfloor\delta\rfloor$ is fixed and the dimension increases, as does the number of steps. For the simulation we set $\aleph=100$, $\varepsilon=0.01$ and the level height $L=5$.\\
\centerline{\includegraphics[scale=0.5]{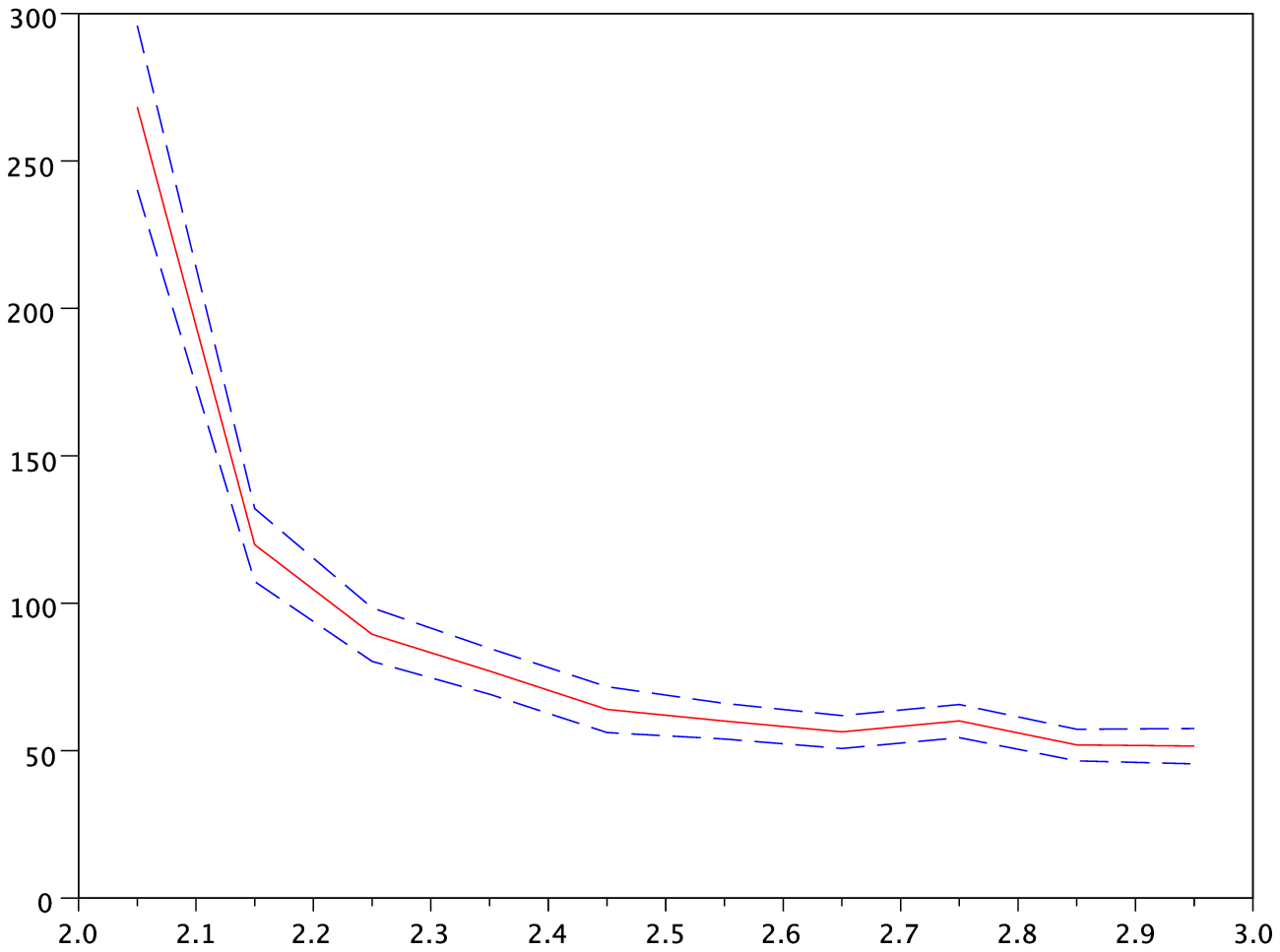}\hspace*{0.4cm}
\includegraphics[scale=0.5]{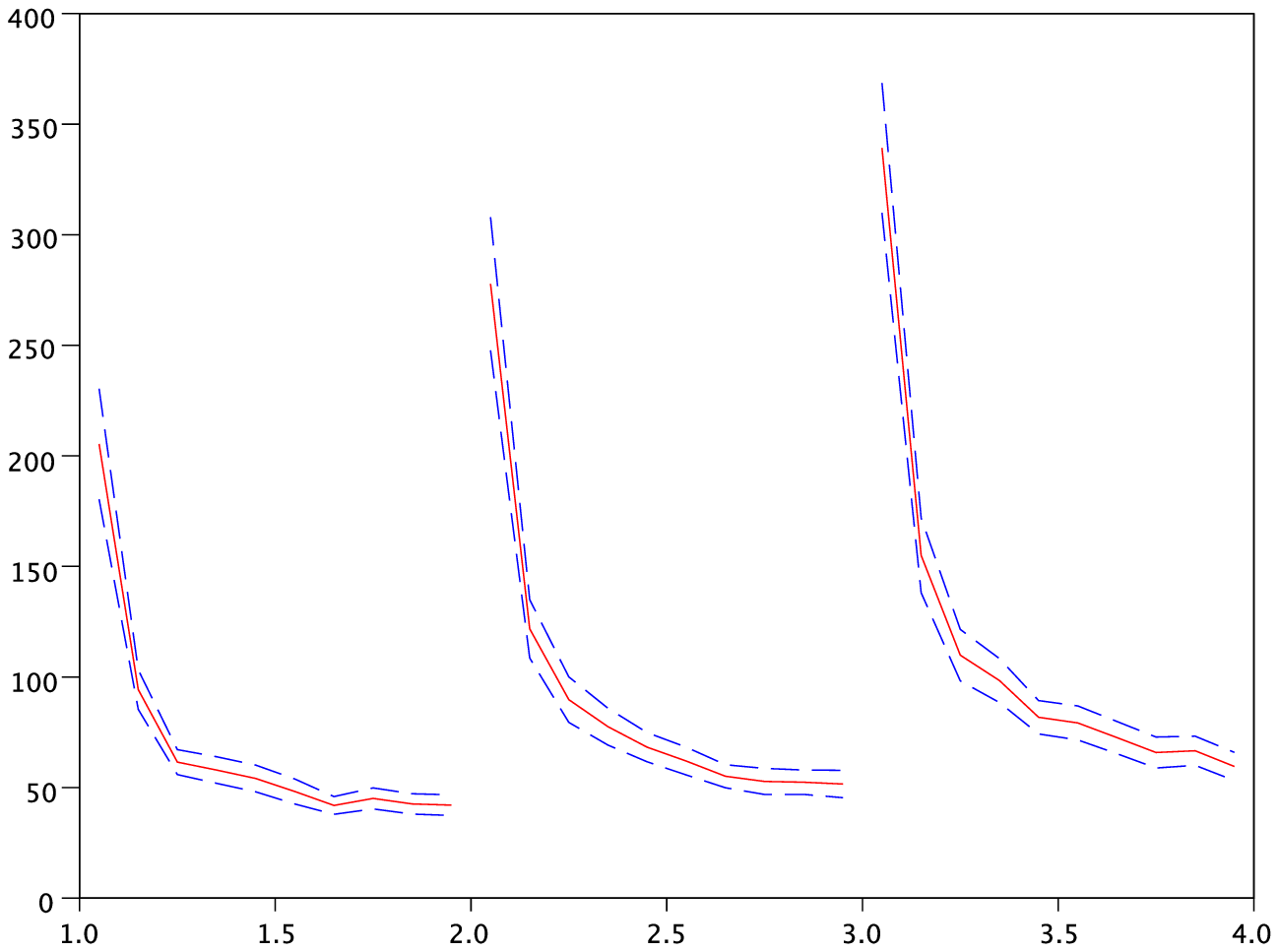}}
\centerline{Average number of step versus $\delta$}
\vspace*{0.2cm}

Now observe the averaged proportion $\mathbb{E}[N^\varepsilon_{\rm integer}/N^\varepsilon]$ as the dimension of the Bessel process increases. It is obvious that this proportion seems to depend mainly on the fractional part of the dimension...\\
 \centerline{\includegraphics[scale=0.5]{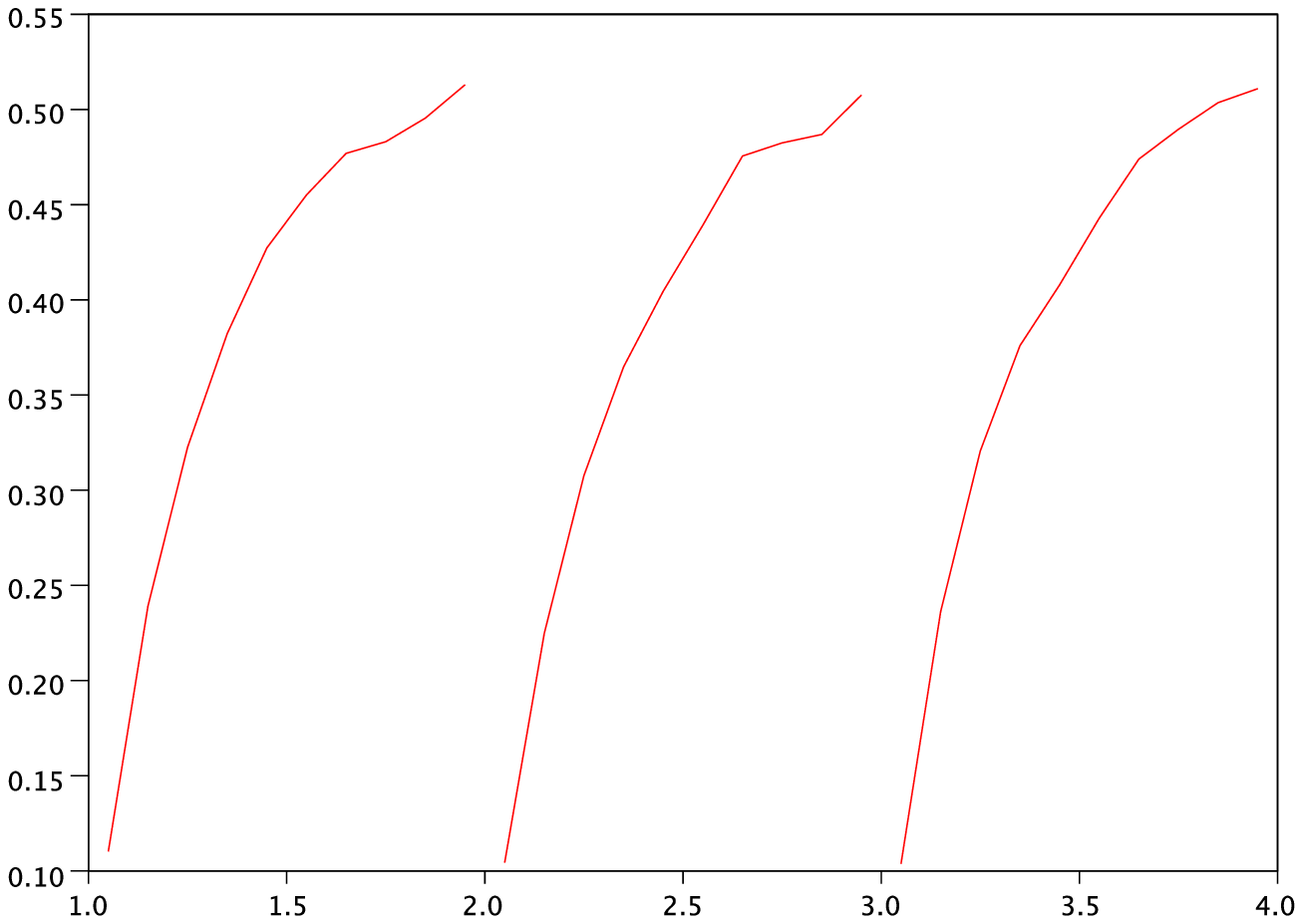}\hspace*{0.4cm}
 \includegraphics[scale=0.5]{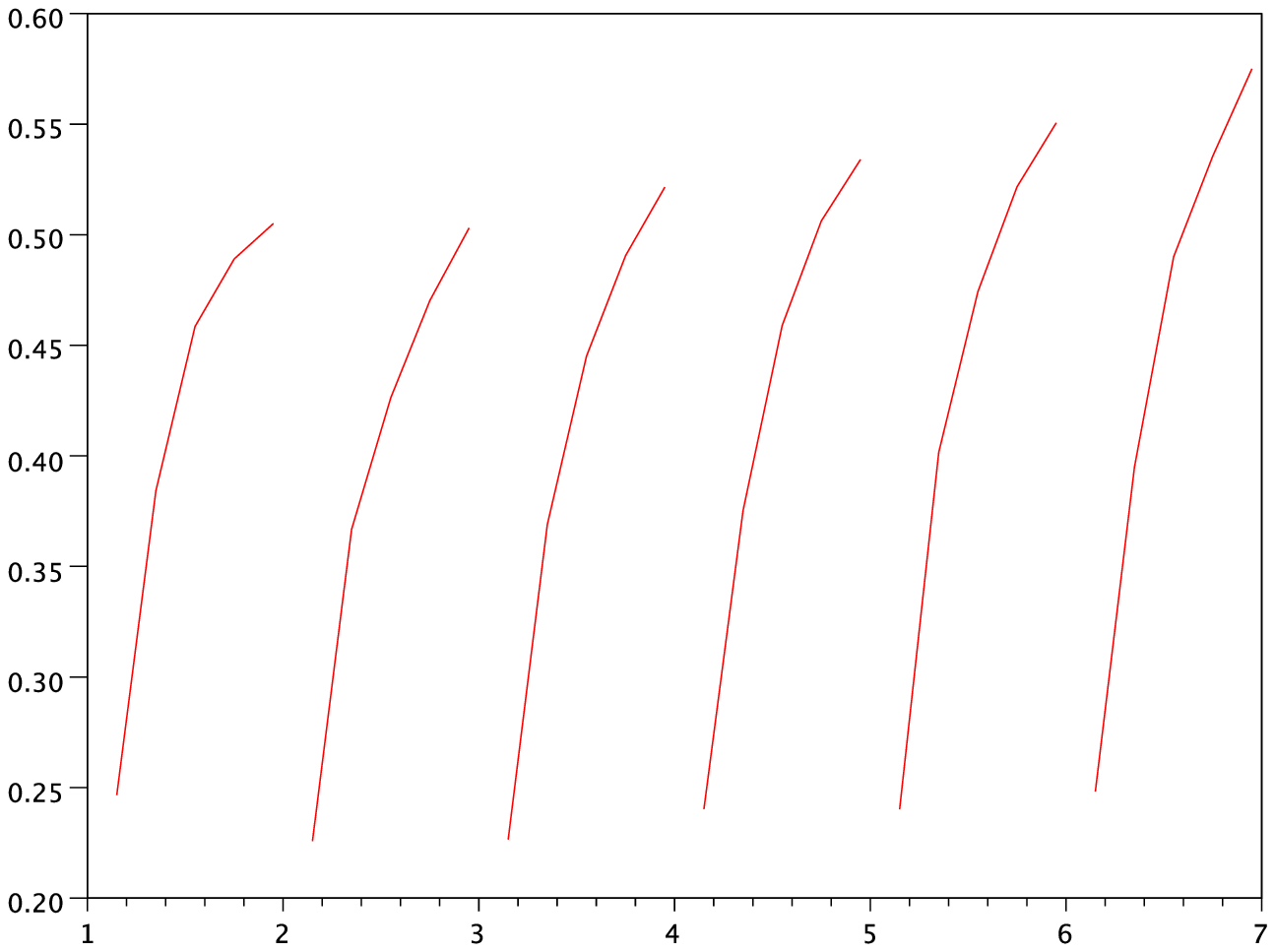}}
 \centerline{Averaged proportion $\mathbb{E}[N^\varepsilon_{\rm integer}/N^\varepsilon]$ versus $\delta$.}\\
\subsection{Number of steps versus the level height}
\par\noindent $\,$

In the previous simulations the level height to reach was always equal to $5$. Let us now present the dependence of the number of steps with respect to this level. The numerical results are obtained for $\aleph=1000$, $\varepsilon=0.01$ and two different Bessel processes (of dimension $\delta=3.8$ and $5.2$). Let us note that this dependence is sub-linear and quite weak, the dimension of the Bessel seems to play a more important role. Observe also that $\varepsilon$ is an upper-bound of $L^2-(M(N^\varepsilon))^2$. We deduce that 
\[
L-M(N^\varepsilon)\le \frac{
\varepsilon}{2L},
\]
and therefore the error of the approximation becomes smaller as $L$ becomes larger. This particular remark is also emphasized by the dependence of $L$ in the bounds \eqref{eq:doubleineq}. Which is why we present a third figure for which $\varepsilon/L$ is fixed and equal to $0.01$.
\\
 \centerline{\includegraphics[scale=0.5]{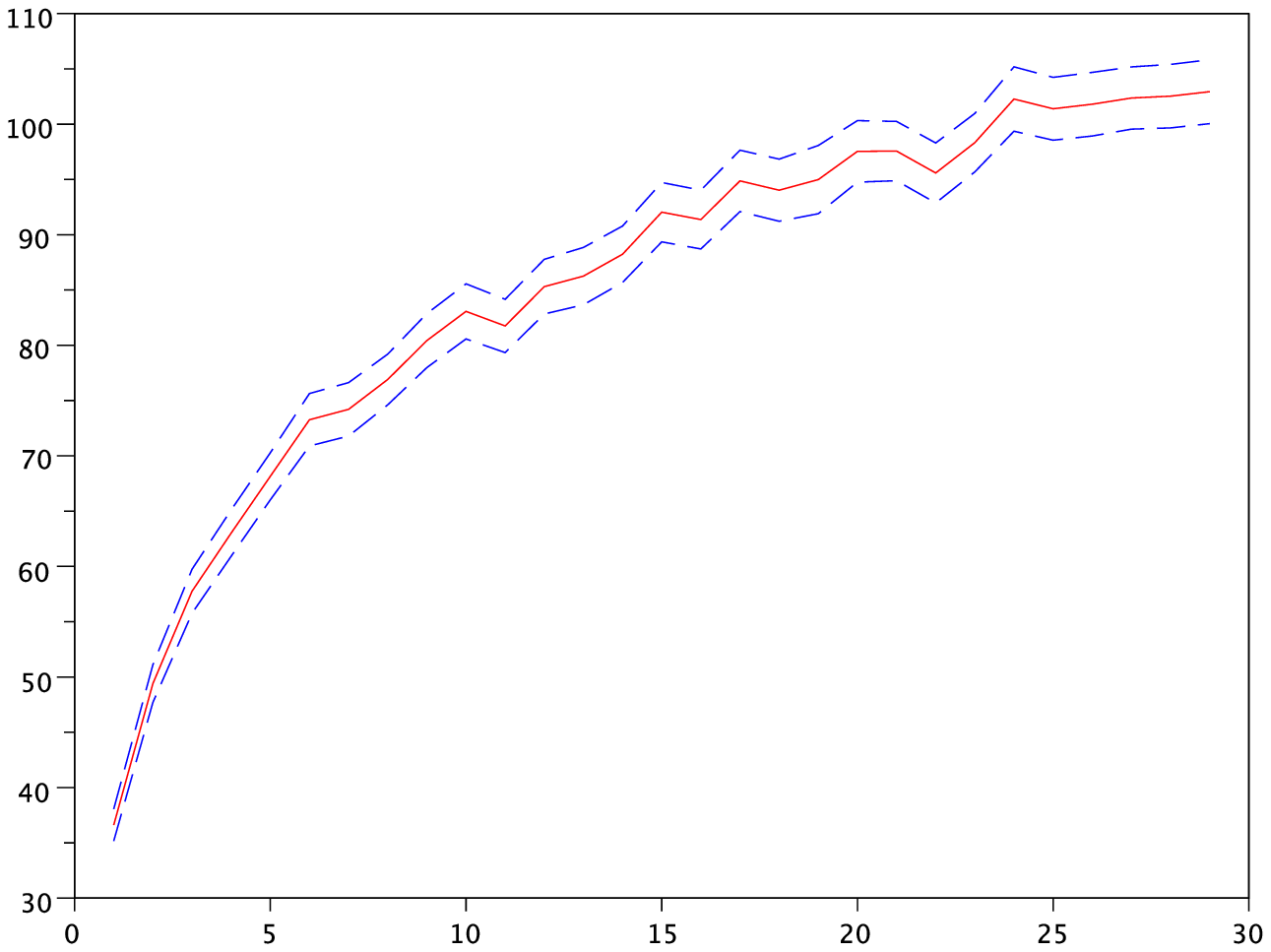}\hspace*{0.4cm}
 \includegraphics[scale=0.5]{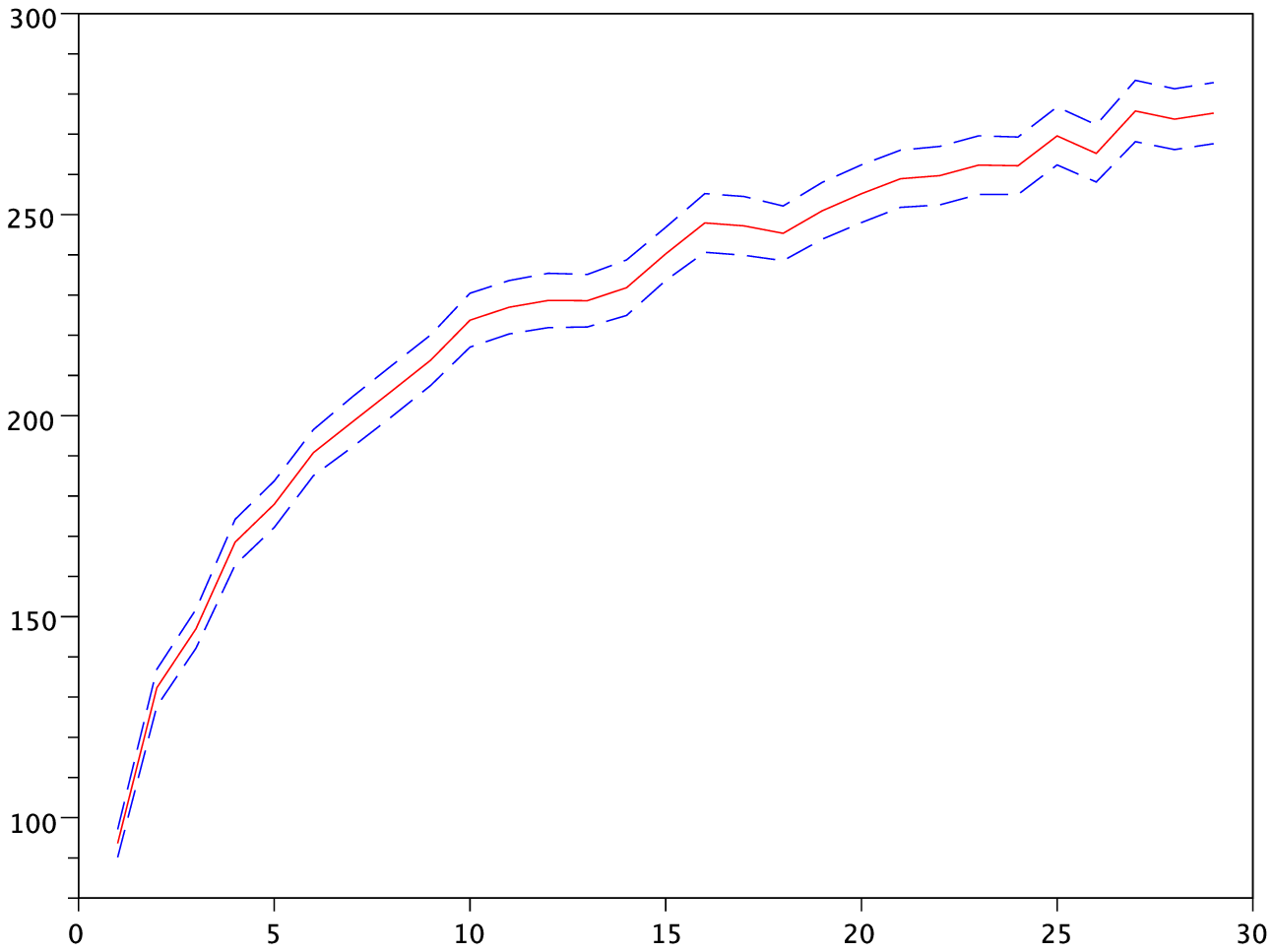}}
 \centerline{Averaged number of step versus the level $L$ (dimension $\delta=3.8$ and $\delta=5.2$)}
 \centerline{\includegraphics[scale=0.5]{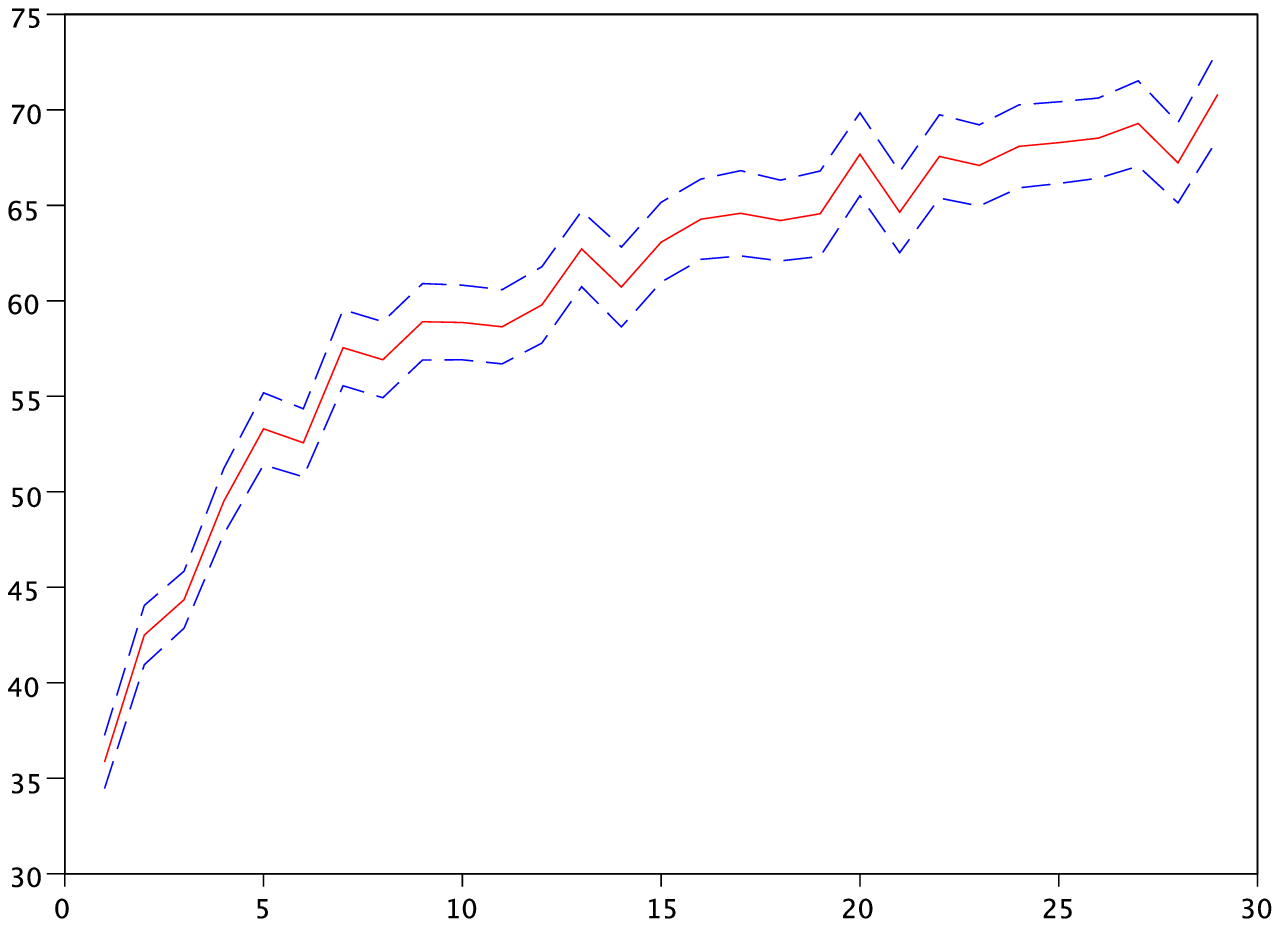}}
 \centerline{Averaged number of steps versus $L$ for fixed $\varepsilon/L$ and $\delta=3.8$}\\
\subsection{Number of generated random variables}
\par\noindent $\,$

Finally let us study the number of random variables used in the simulation of the Bessel hitting times.  Each step of the algorithm (NI) requires a lot of uses of the uniform random generator in order to simulate the first coordinate of the uniform variable on the sphere of dimension $\lfloor \delta\rfloor$, the Gamma distributed variable which appears in the simulation of the hitting times of curved boundaries (here we use Johnk's algorithm, see for instance \cite{Devroye}, page 418), and finally the rejection method for the condition law  described in the \emph{Realization of the algorithm}. The following figures concern simulations with the parameters $L=5$ and $\varepsilon=0.01$ \\
\centerline{\includegraphics[scale=0.5]{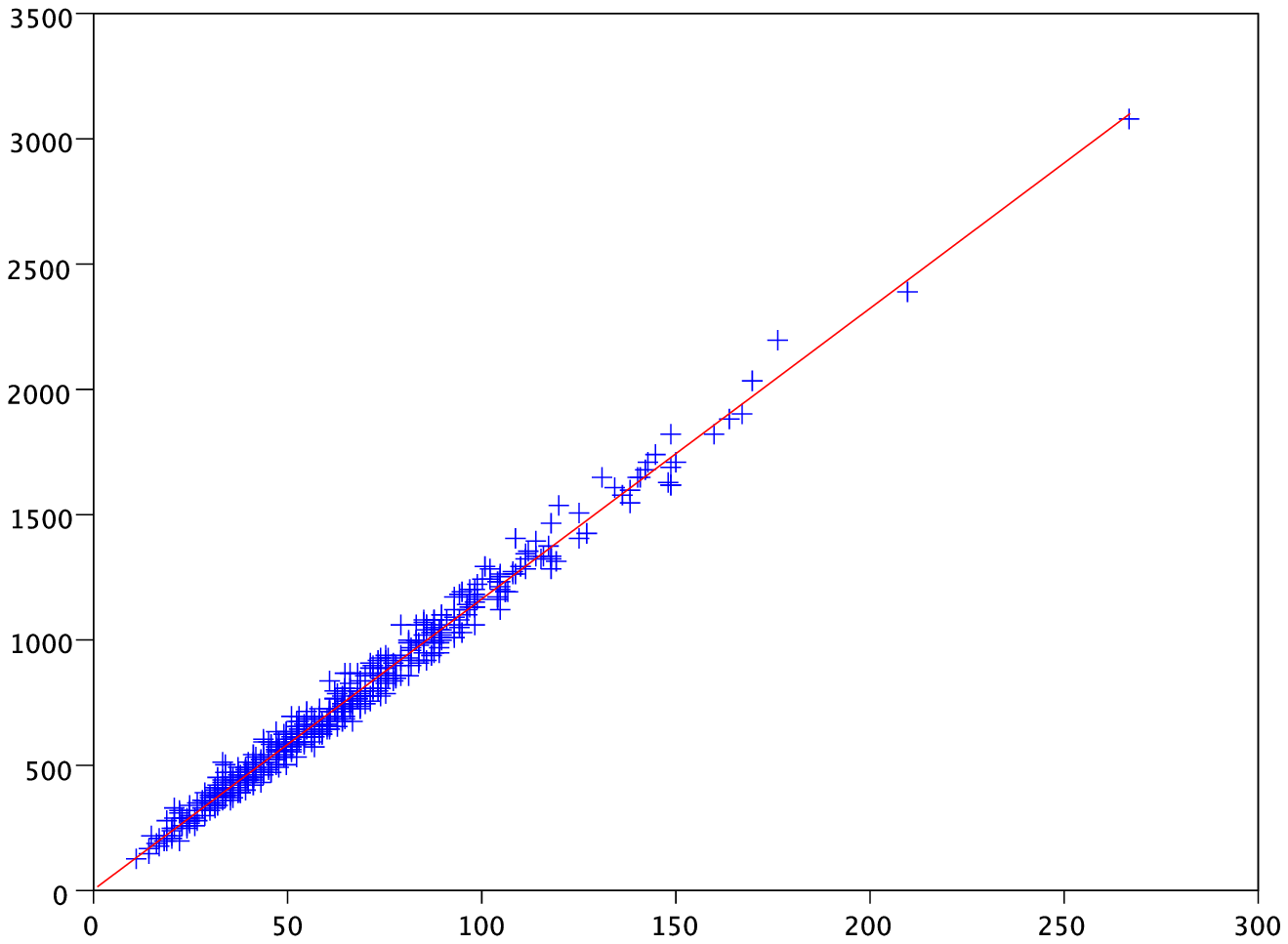}\hspace*{0.4cm}
\includegraphics[scale=0.5]{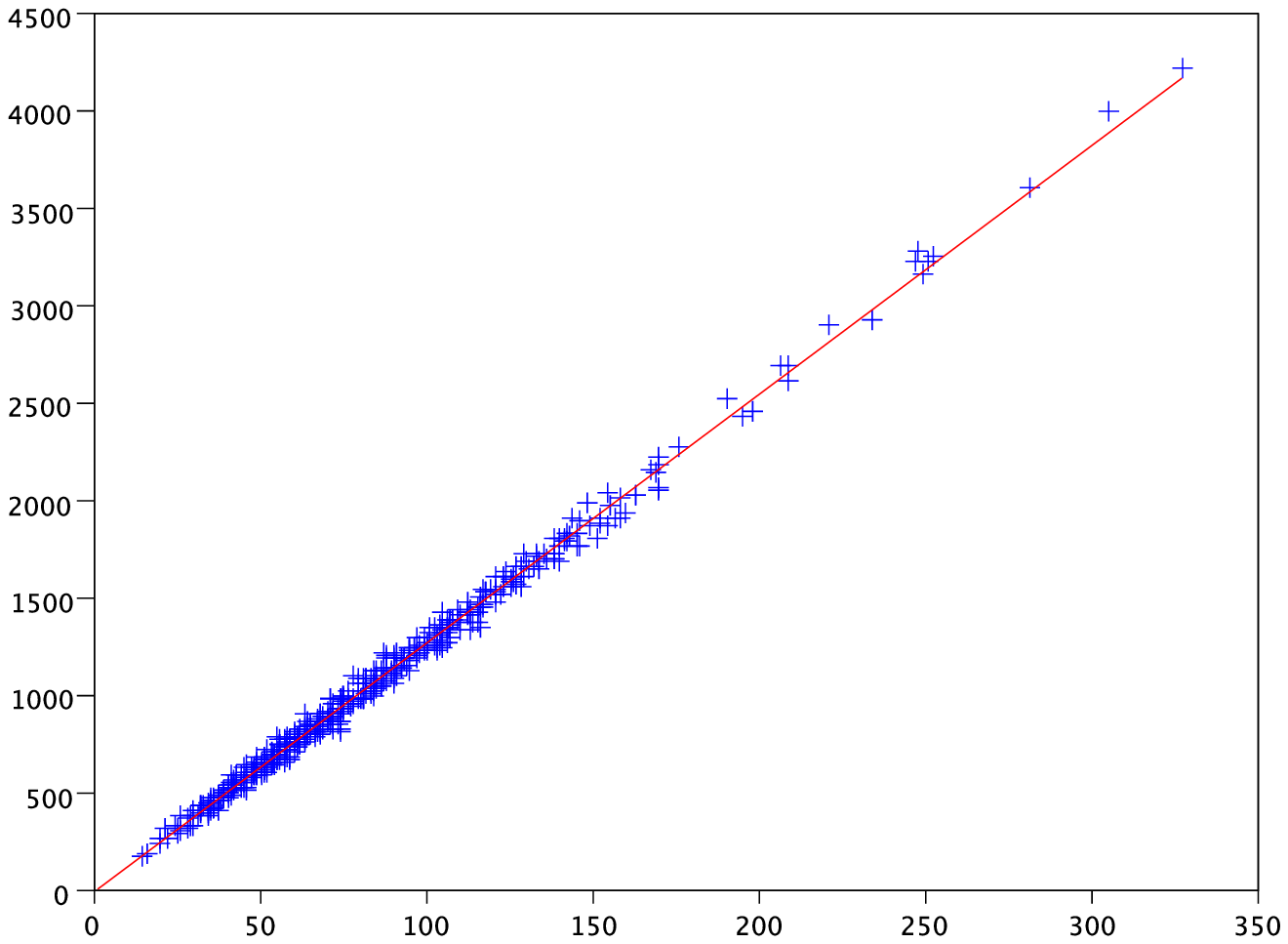}}
 \centerline{Number of random variables required versus number of steps for $\delta=2.5$ and  $\delta=4.8$ }\\

Let us end the numerical section by noting that the algorithm (NI) is quite difficult to use when the fractional part of the dimension $\delta$ i.e. $\delta-\lfloor \delta\rfloor$ is small, the number of steps becomes huge. Moreover, the parameter $\gamma$ appearing in the algorithm (NI) is needed for technical reason and influences the number of steps, so we suggest choosing $\gamma$ as close as possible to $1$.\\ 

\bibliographystyle{plain}
\bibliography{htbiblio1}
\end{document}